\documentclass[a4paper]{amsart}

\usepackage[german,english]{babel}
\usepackage{graphicx}
\RequirePackage{amsmath}

\RequirePackage{amsmath}
\RequirePackage{bm}
\RequirePackage{amssymb}
\RequirePackage{upref}
\RequirePackage{amsthm}
\RequirePackage{enumerate}
\RequirePackage{pb-diagram}
\RequirePackage{amsfonts}
\RequirePackage[mathscr]{eucal}
\RequirePackage{verbatim}
\RequirePackage{xr}
\RequirePackage{graphicx}
\usepackage{calc}
\usepackage{xspace}

\newcommand{\cf}{cf.\@\xspace}
\newcommand{\resp}{resp.\@\xspace}

\newcommand{\aev}{a.e.\@\xspace}

\newcommand{\wlogc}{w.l.o.g.\@\xspace}

\newcommand{\al}{\alpha}
\newcommand{\bet}{\beta}
\newcommand{\ga}{\gamma}
\newcommand{\de}{\delta }
\newcommand{\e}{\epsilon}

\newcommand{\f}{\varphi}
\newcommand{\h}{\eta}

\newcommand{\ka}{\kappa}
\newcommand{\lam}{\lambda}

\newcommand{\n}{\nu}
\newcommand{\om}{\omega}

\newcommand{\s}{\sigma}
\newcommand{\x}{\xi}

\newcommand{\C}{\varGamma}
\newcommand{\D}{\varDelta}
\newcommand{\F}{\varPhi}
\newcommand{\Lam}{\varLambda}
\newcommand{\Om}{\varOmega}

\newcommand{\di}[1]{#1\nobreakdash-\hspace{0pt}dimensional}

\newcommand{\fu}[3]{#1\hspace{0pt}_{|_{#2_#3}}}
\newcommand{\fv}[2]{#1\hspace{0pt}_{|_{#2}}}

\newcommand{\so}{{\mc S_0}}

\newcommand{\const}{\tup{const}}
\newcommand{\ndash}{\nobreakdash--}

\newcommand{\msp[1]}[1]{\mspace{#1mu}}

\newcommand{\R}[1][n+1]{{\protect\mathbb R}^{#1}}

\newcommand{\N}{{\protect\mathbb N}}

\newcommand{\eR}{\stackrel{\lower1ex \hbox{\rule{6.5pt}{0.5pt}}}{\msp[3]\R[]}}
\newcommand{\eN}{\stackrel{\lower1ex \hbox{\rule{6.5pt}{0.5pt}}}{\msp[1]\N}}
\newcommand{\eO}{\stackrel{\lower1ex
\hbox{\rule{6pt}{0.5pt}}}{\msc O}}

\DeclareMathOperator{\diam}{diam}

\DeclareMathOperator{\graph}{graph}

\newcommand\ra{\rightarrow}

\newcommand\pa{\partial}
\newcommand\pde[2]{\frac {\partial#1}{\partial#2}}
 
\newcommand\sql[1][u]{\sqrt{1-|D#1|^2}}
\newcommand{\un}{\infty}
\newcommand{\A}{\forall}

\newcommand{\set}[2]{\{\,#1\colon #2\,\}}
\newcommand{\uu}{\cup}
\newcommand{\ii}{\cap}
\newcommand{\uuu}{\bigcup}

\newcommand{\uud}{ \stackrel{\lower 1ex \hbox {.}}{\uu}}
\newcommand{\uuud}[1]{ \stackrel{\lower 1ex \hbox {.}}{\uuu_{#1}}}
\newcommand\su{\subset}

\newcommand{\sminus}[1][28]{\raise 0.#1ex\hbox{$\scriptstyle\setminus$}}

\newcommand{\wed}{\wedge}

\newcommand\ti{\times }

\newcommand{\abs}[1]{\lvert#1\rvert}

\newcommand{\norm}[1]{\lVert#1\rVert}

\newcommand{\nnorm}[1]{| \mspace{-2mu} |\mspace{-2mu}|#1| \mspace{-2mu}
|\mspace{-2mu}|}
\newcommand{\spd}[2]{\protect\langle #1,#2\protect\rangle}

\newcommand\ch[3]{\varGamma_{#1#2}^#3}
\newcommand\cha[3]{{\bar\varGamma}_{#1#2}^#3}

\newcommand{\riem}[4]{R_{#1#2#3#4}}
\newcommand{\riema}[4]{{\bar R}_{#1#2#3#4}}

\newcommand{\tit}{\textit}

\newcommand{\tup}{\textup}

\newcommand{\mc}{\protect\mathcal}
\newcommand{\msc}{\protect\mathscr}

\providecommand{\bysame}{\makebox[3em]{\hrulefill}\thinspace}

\newcommand{\cq}[1]{\glqq{#1}\grqq\,}
\newcommand{\cqp}[1]{\glqq{\ignorespaces #1\ignorespaces}\grqq}

\newcommand{\bt}{\begin{thm}}
\newcommand{\bl}{\begin{lem}}
\newcommand{\bc}{\begin{cor}}
\newcommand{\bd}{\begin{definition}}
\newcommand{\bpp}{\begin{prop}}
\newcommand{\br}{\begin{rem}}
\newcommand{\bn}{\begin{note}}
\newcommand{\be}{\begin{ex}}
\newcommand{\bes}{\begin{exs}}
\newcommand{\bb}{\begin{example}}
\newcommand{\bbs}{\begin{examples}}
\newcommand{\ba}{\begin{axiom}}
\newcommand{\bas}{\begin{assumption}}

\newcommand{\et}{\end{thm}}
\newcommand{\el}{\end{lem}}
\newcommand{\ec}{\end{cor}}
\newcommand{\ed}{\end{definition}}
\newcommand{\epp}{\end{prop}}
\newcommand{\er}{\end{rem}}
\newcommand{\en}{\end{note}}
\newcommand{\ee}{\end{ex}}
\newcommand{\ees}{\end{exs}}
\newcommand{\eb}{\end{example}}
\newcommand{\ebs}{\end{examples}}
\newcommand{\ea}{\end{axiom}}
\newcommand{\eas}{\end{assumption}}

\newcommand{\bp}{\begin{proof}}
\newcommand{\ep}{\end{proof}}
\newcommand{\eps}{\renewcommand{\qed}{}\end{proof}}

\newcommand{\bal}{\begin{align}}

\newcommand{\bi}[1][1.]{\begin{enumerate}[\upshape #1]}
\newcommand{\bia}[1][(1)]{\begin{enumerate}[\upshape #1]}
\newcommand{\bin}[1][1]{\begin{enumerate}[\upshape\bfseries #1]}
\newcommand{\bir}[1][(i)]{\begin{enumerate}[\upshape #1]}
\newcommand{\bic}[1][(i)]{\begin{enumerate}[\upshape\hspace{2\cma}#1]}
\newcommand{\bis}[2][1.]{\begin{enumerate}[\upshape\hspace{#2\parindent}#1]}
\newcommand{\ei}{\end{enumerate}}

\newcommand\ndots{\raise 0.47ex \hbox {,}\hskip0.06em\cdots %
     \raise 0.47ex \hbox {,}\hskip0.06em} 

\newcommand{\q}{\quad}
\newcommand{\qq}{\qquad}

\newcommand{\hp}{\hphantom}

\newcommand\nd{\noindent}

\newskip\Csmallskipamount                                                
\Csmallskipamount=\smallskipamount
\newskip\Cmedskipamount
\Cmedskipamount=\medskipamount
\newskip\Cbigskipamount
\Cbigskipamount=\bigskipamount

\newcommand\cvs{\vspace\Csmallskipamount}   
\newcommand\cvm{\vspace\Cmedskipamount}

\newskip\csa
\csa=\smallskipamount

\newskip\cma
\cma=\medskipamount

\newskip\cba
\cba=\bigskipamount

\newdimen\spt
\newcommand\citem{\cvs\advance\itemno by
1{(\romannumeral\the\itemno})\hskip3pt}
\newcommand{\bitem}{\cvm\nd\advance\itemno by
1{\bf\the\itemno}\hspace{\cma}}

\newcount\itemno
\newcommand{\las}[1]{\label{S:#1}}

\newcommand{\lae}[1]{\label{E:#1}}
\newcommand{\lat}[1]{\label{T:#1}}
\newcommand{\lal}[1]{\label{L:#1}}
\newcommand{\lad}[1]{\label{D:#1}}

\newcommand{\lan}[1]{\label{N:#1}}

\newcommand{\lar}[1]{\label{R:#1}}

\newcommand{\rs}[1]{Section~\ref{S:#1}}

\newcommand{\rt}[1]{Theorem~\ref{T:#1}}
\newcommand{\rl}[1]{Lemma~\ref{L:#1}}
\newcommand{\rd}[1]{Definition~\ref{D:#1}}

\newcommand{\re}[1]{\eqref{E:#1}}
\newcommand{\frt}[1]{Theorem~\ref{T:#1} on page~\tup{\pageref{T:#1}}}
\newcommand{\frl}[1]{Lemma~\ref{L:#1} on page~\tup{\pageref{L:#1}}}
\newcommand{\frr}[1]{Remark~\ref{R:#1} on page~\tup{\pageref{R:#1}}}
\newcommand{\frd}[1]{Definition~\ref{D:#1} on page~\tup{\pageref{D:#1}}}

\newskip\thmskip
\thmskip=\parindent

\newskip\hsk
\newenvironment{hinw}{\labelsep=0pt\begin{list}{}{\labelsep=0pt\itemindent=0pt\labelwidth=0pt\leftmargin=\parindent\rightmargin=0pt\partopsep=\cba}%
\item\it\nopagebreak\nopagebreak}%
{\end{list}}

\newcommand\bh{\begin{hinw}}
\newcommand{\eh}{\end{hinw}}

\newtheoremstyle{normal}
  {\cba}
  {\cba}
  {}
  {\thmskip}
  {\bfseries}
  {.}
  {\hsk}
  {}

\newtheoremstyle{abschnitt}
  {\cba}
  {\cba}
  {}
  {\thmskip}
  {\bfseries}
  {.}
  {\hsk}
  {}

\newtheoremstyle{italic}
  {\cba}
  {\cba}
  {\itshape}
  {\thmskip}
  {\bfseries}
  {.}
  {\hsk}
  {}

\newtheoremstyle{aufgaben}
  {\cba}
  {\cba}
  {}
  {}
  {\normalsize\bfseries}
  {.}
  {\hsk}
  {}

\newtheoremstyle{break}
  {\cba}
  {\cba}
  {\itshape}
  {}
  {\bfseries}
  {.}
  {\newline}
  {}

\swapnumbers
\theoremstyle{italic}
\newtheorem{thm}[subsection]{Theorem}
\newtheorem{lem}[subsection]{Lemma}
\newtheorem{prop}[subsection]{Proposition}
\newtheorem{cor}[subsection]{Corollary}

\theoremstyle{normal}
\newtheorem{rem}[subsection]{Remark}
\newtheorem{definition}[subsection]{Definition}
\newtheorem{example}[subsection]{Example}
\newtheorem{examples}[subsection]{Examples}
\newtheorem{ex}[subsection]{Exercise}
\newtheorem{note}[subsection]{}
\newtheorem{axiom}[subsection]{Axiom}
\newtheorem{assumption}[subsection]{Assumption}

\theoremstyle{aufgaben}
\newtheorem{exs}[subsection]{Exercises}

\swapnumbers

\numberwithin{equation}{section}
\numberwithin{figure}{section}

\newenvironment{textequation}[1][0.8]
{\begin{equation}
\begin{aligned}
\begin{minipage}{#1\linewidth}}
{\end{minipage}
\end{aligned}
\end{equation}
\ignorespacesafterend}

\newcommand{\btext}{\begin{textequation}}
\newcommand{\etext}{\end{textequation}}
\makeatletter
\def\hinweis{\@startsection{subsection}{2}%
 \z@{0.7\linespacing\@plus 0.5\linespacing}{0.7\linespacing}%
{\normalfont\itshape\indent}}
\makeatother
\RequirePackage{bm}
\RequirePackage{amssymb}
\RequirePackage{upref}
\RequirePackage{amsthm}
\RequirePackage{enumerate}
\RequirePackage{pb-diagram}
\RequirePackage{amsfonts}
\RequirePackage[mathscr]{eucal}
\RequirePackage{verbatim}
\RequirePackage{xr}
\RequirePackage{graphicx}
\usepackage{calc}
\usepackage{xspace}

\newcommand{\ind}[1]{#1}
\newcommand{\indexs}[1]{\relax}
\renewcommand{\index}[1]{\relax}
\newcommand{\inds}[1]{#1}
\newcommand{\fre}[1]{\eqref{E:#1} on page~\tup{\tup{\pageref{E:#1}}}}
\newcommand{\lab}[1]{\label{B:#1}}

\newcommand{\frs}[1]{Section~\ref{S:#1} on page~\tup{\pageref{S:#1}}}

\RequirePackage{upref}
\RequirePackage{amsthm}
\RequirePackage{enumerate}
\usepackage[mathscr]{eucal}

\usepackage{xr-hyper}

\usepackage{calc}

\newlength{\oddsidemarginlength}
\newlength{\topmarginlength}

\newcounter{numberoflines}
\newcounter{tempcc}
\oddsidemargin=\oddsidemarginlength
\evensidemargin=\oddsidemargin
\topmargin=\topmarginlength
\begin{document}

\flushbottom


\title{Curvature flows in semi-Riemannian manifolds}

\author{Claus Gerhardt}
\address{Ruprecht-Karls-Universit\"at, Institut f\"ur Angewandte Mathematik,
Im Neuenheimer Feld 294, 69120 Heidelberg, Germany}
\email{gerhardt@math.uni-heidelberg.de}
\urladdr{http://www.math.uni-heidelberg.de/studinfo/gerhardt/}
\thanks{This research was supported by the Deutsche Forschungsgemeinschaft.}

%
\subjclass[2000]{35J60, 53C21, 53C44, 53C50, 58J05}
\keywords{semi-Riemannian manifold, mass, stable solutions, cosmological spacetime, general relativity, curvature flows, ARW spacetime}
\date{\today}
%


\begin{abstract}
We prove that the limit hypersurfaces of converging curvature flows are stable, if the initial velocity has a weak sign, and give a survey of the existence and regularity results.
\end{abstract}

\maketitle

\tableofcontents

\setcounter{section}{0}
\section{Introduction}
In this paper we want to give a survey of the existence and regularity results for extrinsic curvature flows in semi-Riemannian manifolds, i.e., Riemannian or Lorentzian ambient spaces, with an emphasis on flows in Lorentzian spaces. In order to treat both cases simultaneously terminology like spacelike, timelike, etc., that only makes sense in a Lorentzian setting should be ignored in the Riemannian case.

The general stability result for the limit hypersurfaces of converging curvature flows in \rs{5c} is new. The regularity result in \rt{6.5.5}---especially the time independent $C^{m+2,\al}$-estimates---for converging curvature flows that are graphs  is interesting too.

\section{Notations and preliminary results}\las 2

The main objective of this section is to state the equations of Gau{\ss}, Codazzi,
and Weingarten for hypersurfaces. In view of the subtle but important
difference  that is to be seen in the \tit{Gau{\ss} equation} depending on the
nature of the ambient space---Riemannian or Lorentzian---, which we already
mentioned in the introduction, we shall formulate the governing equations of a
hypersurface $M$ in a semi-Riemannian \di{(n+1)} space $N$, which is either
Riemannian or Lorentzian. Geometric quantities in $N$ will be denoted by
$(\bar g_{\alpha\beta}),(\riema \alpha\beta\gamma\delta)$, etc., and those in $M$ by $(g_{ij}), (\riem
ijkl)$, etc. Greek indices range from $0$ to $n$ and Latin from $1$ to $n$; the
summation convention is always used. Generic coordinate systems in $N$ resp.
$M$ will be denoted by $(x^\alpha)$ resp. $(\x^i)$. Covariant differentiation will
simply be indicated by indices, only in case of possible ambiguity they will be
preceded by a semicolon, i.e. for a function $u$ in $N$, $(u_\alpha)$ will be the
gradient and
$(u_{\alpha\beta})$ the Hessian, but e.g., the covariant derivative of the curvature
tensor will be abbreviated by $\riema \alpha\beta\gamma{\delta;\e}$. We also point out that
\begin{equation}
\riema \alpha\beta\gamma{\delta;i}=\riema \alpha\beta\gamma{\delta;\e}x_i^\e
\end{equation}
with obvious generalizations to other quantities.

Let $M$ be a \tit{spacelike} hypersurface, i.e. the induced metric is Riemannian,
with a differentiable normal $\n$. We define the signature of $\n$, $\sigma=\s(\n)$, by
\begin{equation}
\sigma=\bar g_{\alpha\beta}\n^\alpha\n^\beta=\spd \n\n.
\end{equation}
In case $N$ is Lorentzian, $\sigma=-1$, and $\n$ is timelike.

In local coordinates, $(x^\alpha)$ and $(\x^i)$, the geometric quantities of the
spacelike hypersurface $M$ are connected through the following equations
\begin{equation}\lae{2.3}
x_{ij}^\alpha=-\sigma h_{ij}\n^\alpha
\end{equation}
the so-called \tit{Gau{\ss} formula}. Here, and also in the sequel, a covariant
derivative is always a \tit{full} tensor, i.e.

\begin{equation}
x_{ij}^\alpha=x_{,ij}^\alpha-\ch ijk x_k^\alpha+\cha \beta\gamma\alpha x_i^\beta x_j^\gamma.
\end{equation}
The comma indicates ordinary partial derivatives.

In this implicit definition the \tit{second fundamental form} $(h_{ij})$ is taken
with respect to $-\sigma\n$.

The second equation is the \tit{Weingarten equation}
\begin{equation}
\n_i^\alpha=h_i^k x_k^\alpha,
\end{equation}
where we remember that $\n_i^\alpha$ is a full tensor.

Finally, we have the \tit{Codazzi equation}
\begin{equation}
h_{ij;k}-h_{ik;j}=\riema\alpha\beta\gamma\delta\n^\alpha x_i^\beta x_j^\gamma x_k^\delta
\end{equation}
and the \tit{Gau{\ss} equation}
\begin{equation}
\riem ijkl=\sigma \{h_{ik}h_{jl}-h_{il}h_{jk}\} + \riema \alpha\beta\gamma\delta x_i^\alpha x_j^\beta x_k^\gamma
x_l^\delta.
\end{equation}
Here, the signature of $\n$ comes into play.

\bd\lad{1.6.7}
(i) Let $F\in C^0(\bar\C)\ii C^{2,\al}(\C)$ be a strictly monotone curvature function, where $\C\su\R[n]$ is a convex, open, symmetric cone containing the positive cone, such that
\begin{equation}
\fv F{\pa \C}=0\q\wed\q \fv F\C>0.
\end{equation}

Let $N$ be semi-Riemannian. A spacelike, orientable\footnote{A hypersurface is said to be orientable, if it has a continuous normal field.} hypersurface $M\su N$ is called \tit{admissible}\index{admissible hypersurface}, if its principal curvatures with respect to a chosen normal lie in $\C$. This definition also applies to subsets of $M$.

\cvm
(ii) Let $M$ be an admissible hypersurface and $f$ a function defined in a neighbourhood of $M$. $M$ is said to be an \tit{\ind{upper barrier}} for the pair $(F,f)$, if
\begin{equation}
\fv FM\ge f
\end{equation}

(iii) Similarly, a spacelike, orientable hypersurface $M$ is called a \tit{\ind{lower barrier}} for the pair $(F,f)$, if at the points $\varSigma\su M$, where
$M$ is admissible, there holds
\begin{equation}
\fv F\varSigma\le f.
\end{equation}
$\varSigma$ may be empty.

\cvm
(iv) If we consider the mean curvature function, $F=H$, then we suppose $F$ to be defined in $\R[n]$ and any spacelike, orientable hypersurface is admissible.
\ed

One of the assumptions that are used when proving a priori estimates is that there exists a strictly convex function 
$\chi\in C^2(\bar \Omega)$ in a given domain $\Om$. We shall state sufficient geometric conditions
guaranteeing the existence of such a function. The lemma below will be valid in Lorentzian as well as Riemannian manifolds, but we formulate and prove it only for the Lorentzian case.

\bl\lal{2.8.1}
Let $N$ be globally hyperbolic, $\so$ a Cauchy hypersurface, $(x^\alpha)$ a special
coordinate system associated with $\so$, and $\bar \Omega\su N$ be compact. Then,
there exists a strictly convex function $\chi\in C^2(\bar \Omega)$ provided the level
hypersurfaces $\{x^0=\const\}$ that intersect $\bar \Omega$ are strictly convex.
\el

\bp
For greater clarity set $t=x^0$, i.e., $t$ is a globally defined time function. Let
$x=x(\x)$ be a local representation for $\{t=\const\}$, and $t_i,t_{ij}$ be the
covariant derivatives of $t$ with respect to the induced metric, and
$t_\alpha,t_{\alpha\beta}$ be the covariant derivatives in $N$, then
\begin{equation}
0=t_{ij}=t_{\alpha\beta}x_i^\alpha x_j^\beta+t_\alpha x_{ij}^\alpha,
\end{equation}
and therefore,
\begin{equation}\lae{2.26}
t_{\alpha\beta}x_i^\alpha x_j^\beta=-t_\alpha x_{ij}^\alpha=-\bar h_{ij} t_\alpha \n^\alpha.
\end{equation}
Here, $(\n^\alpha)$ is past directed, i.e., the right-hand side in \re{2.26} is positive
definite in $\bar \Omega$, since $(t_\alpha)$ is also past directed.

Choose $\lambda>0$ and define $\chi=e^{\lambda t}$, so that
\begin{equation}
\chi_{\alpha\beta}=\lambda^2e^{\lambda t} t_\alpha t_\beta+\lambda e^{\lambda t} t_{\alpha\beta}.
\end{equation}

Let $p\in \Omega$ be arbitrary, $\mc S=\{t=t(p)\}$ be the level hypersurface through
$p$, and $(\h^\alpha)\in T_p(N)$. Then, we conclude
\begin{equation}
e^{-\lambda t}\chi_{\alpha\beta}\h^\alpha \h^\beta=\lambda^2\abs{\h^0}^2+\lambda t_{ij}\h^i\h^j+2\lambda t_{0j}\h^0\h^i,
\end{equation}
where $t_{ij}$ now represents the left-hand side in  \re{2.26}, and we infer further
\begin{equation}
\begin{aligned}
e^{-\lambda t}\chi_{\alpha\beta}\h^\alpha\h^\beta&\ge \tfrac{1}{2} \lambda^2 {\abs\h^0}^2 +[\lambda
\e-c_\e]\sigma_{ij}\h^i\h^j\\
&\ge \tfrac{\e}{2}\lambda \{-\abs{\h^0}^2+\sigma_{ij}\h^i\h^j\}
\end{aligned}
\end{equation}
for some $\e>0$, and where $\lambda$ is supposed to be large. Therefore, we
have in $\bar \Omega$
\begin{equation}
\chi_{\alpha\beta}\ge c\bar g_{\alpha\beta}\,,\q c>0,
\end{equation}
i.e., $\chi$ is strictly convex.
\ep

\section{Evolution equations for some geometric quantities}\las{1.3}

Curvature flows are used for different purposes, they can be merely vehicles to approximate a stationary solution, in which case the flow is driven not only by a curvature function but also by the corresponding right-hand side, an external force, if you like, or the flow is a pure curvature flow driven only by a curvature function, and it is used to analyze the topology of the initial hypersurface, if the ambient space is Riemannian, or the singularities of the ambient space, in the Lorentzian case.

In this section we are treating very general curvature flows\footnote{We emphasize that we are only considering flows driven by the extrinsic curvature not by the intrinsic curvature.} in a semi-Riemann\-ian manifold $N=N^{n+1}$, though we only have the Riemannian or Lorentzian case in mind, such that the flow can be either a pure curvature flow or may also be driven by an external force. The nature of the ambient space, i.e., the signature of its metric, is expressed by a parameter $\s=\pm 1$, such that $\s=1$ corresponds to the Riemannian and $\s=-1$ the Lorentzian case. The parameter $\s$ can also be viewed as the signature  of the normal of the spacelike hypersurfaces, namely, 
\begin{equation}
\s=\spd\n\nu.
\end{equation}

 Properties like spacelike, achronal, etc., however, 
only make sense, when $N$ is Lorentzian and should be ignored otherwise.
 
 We consider a strictly monotone, symmetric, and concave curvature $F\in C^{4,\al}(\C)$, homogeneous of degree $1$, a function $0<f\in C^{4,\al}(\Om)$,  where $\Om\su N$ is an open set, and a real function $\F\in C^{4,\al}(\R[]_+)$ satisfying
\begin{equation}
\dot\F>0\q \tup{and}\q \ddot\F\le 0.
\end{equation}

For notational reasons, let us abbreviate
\begin{equation}
\tilde f=\F(f).
\end{equation}

Important examples of functions $\F$ are
 \begin{equation}
\F(r)=r,\q
\F(r)=\log r,\q \F(r)=-r^{-1}
\end{equation}
or
\begin{equation}
\F(r)=r^{\frac1k},\q \F(r)=-r^{-\frac1k},\qq k\ge1.
\end{equation}
\br
The latter choices are necessary, if the curvature function $F$ is not homogeneous of degree $1$ but of degree $k$, like the symmetric polynomials $H_k$. In this case we would sometimes like to define $F=H_k$ and not $H_k^{1/k}$, since 
\begin{equation}
F^{ij}=\frac{\pa F}{\pa h_{ij}}
\end{equation}
is then divergence free, if the ambient space is a spaceform, \cf \frl{5.8.2}, though on the other hand we need a concave operator for technical reasons, hence we have to take the $k$-th root.
\er
The curvature flow is given by  the evolution problem
\begin{equation}\lae{1.3.5}
\begin{aligned}
\dot x&=-\s(\F-\tilde f)\n,\\
x(0)&=x_0,
\end{aligned}
\end{equation}
where $x_0$ is an embedding of an initial  compact, spacelike
hypersurface $M_0\su \Om$ of class $C^{6,\al}$, $\F=\F(F)$, and $F$ is evaluated at the principal curvatures
of the flow hypersurfaces $M(t)$, or, equivalently, we may assume that $F$
depends on the second fundamental form $(h_{ij})$ and the metric $(g_{ij})$ of
$M(t)$; $x(t)$ is the embedding of $M(t)$ and $\sigma$ the signature of the  normal $\n=\n(t)$, which is identical to the normal used in the \ind{Gaussian formula} \fre{2.3}.

The initial hypersurface should be \tit{admissible}\index{admissible hypersurface}, i.e., its principal curvatures should belong to the convex, symmetric cone $\C\su \R[n]$.

This is a parabolic problem, so short-time existence is guaranteed, \cf \cite[Chapter 2.5]{cg:cp}

There will be a slight ambiguity in the terminology, since we shall call the
evolution parameter \tit{time}, but this lapse shouldn't cause any
misunderstandings, if the ambient space is Lorentzian.

At the moment we consider a sufficiently smooth solution of the initial value problem \re{1.3.5} and want to show how the metric, the second fundamental form, and the
normal vector of the hypersurfaces $M(t)$ evolve. All time derivatives are
\tit{total} derivatives, i.e., covariant derivatives of tensor fields defined over the curve $x(t)$, \cf \cite[Chapter 11.5]{cg:aII}; $t$ is the flow parameter, also referred to as \tit{\ind{time}}, and $(\xi^i)$ are local coordinates of the initial embedding $x_0=x_0(\xi)$ which will also serve as coordinates for the the flow hypersurfaces $M(t)$. The coordinates in $N$ will be labelled $(x^\al)$, $0\le \al\le n$.
\bl[Evolution of the metric]
The metric $g_{ij}$ of $M(t)$ satisfies the evolution equation
\begin{equation}\lae{1.3.6}
\dot g_{ij}=-2\s(\F-\tilde f)h_{ij}.
\end{equation}
\el
\bp
Differentiating
\begin{equation}
g_{ij}=\spd{x_i}{x_j}
\end{equation}
covariantly with respect to $t$ yields
\begin{equation}
\begin{aligned}
\dot g_{ij}&=\spd{\dot x_i}{x_j}+\spd{x_i}{\dot x_j}\\
&=-2 \s (\F-\tilde f)\spd{x_i}{\nu_j}=-2\s (\F-\tilde f)h_{ij},
\end{aligned}
\end{equation}
in view of the Codazzi equations.
\ep

\bl[Evolution of the normal]
The normal vector evolves according to
\begin{equation}\lae{1.3.9}
\dot \n=\nabla_M(\F-\tilde f)=g^{ij}(\F-\tilde f)_i x_j.
\end{equation}
\el

\bp
Since $\nu$ is unit normal vector we have $\dot\nu\in T(M)$. Furthermore, differentiating 
\begin{equation}
0=\spd\nu{x_i}
\end{equation}
with respect to $t$, we deduce
\begin{equation}
\spd{\dot\nu}{x_i}=-\spd{\nu}{\dot x_i}=(\F-\tilde f)_i.\qedhere
\end{equation}
\ep

\bl[Evolution of the second fundamental form]\lal{1.3.3}
The second fundamental form evolves according to
\begin{equation}\lae{1.3.12}
\dot h_i^j=(\F-\tilde f)_i^j+\sigma (\F-\tilde f) h_i^k h_k^j+\sigma (\F-\tilde f) \riema
\alpha\beta\gamma\delta\n^\alpha x_i^\beta \n^\gamma x_k^\delta g^{kj}
\end{equation}
and
\begin{equation}\lae{1.3.13}
\dot h_{ij}=(\F-\tilde f)_{ij}-\sigma (\F-\tilde f) h_i^k h_{kj}+\sigma (\F-\tilde f) \riema
\alpha\beta\gamma\delta\n^\alpha x_i^\beta \n^\gamma x_j^\delta.
\end{equation}
\el

\bp
We use the \tit{\ind{Ricci identities}} to interchange the covariant derivatives of $\nu$ with respect to $t$ and $\xi^i$
\begin{equation}
\begin{aligned}
\tfrac D{dt}(\n^\al_i)&=(\dot\nu^\al)_i-\bar R^\al_{\hp{\al}\bet\ga\de}\nu^\bet x^\ga_i \dot x^\de\\
&= g^{kl}(\F-\tilde f)_{ki}x^\al_l+g^{kl}(\F-\tilde f)_kx^\al_{li}-\bar R^\al_{\hp{\al}\bet\ga\de} \nu^\bet x^\ga_i \dot x^\de.
\end{aligned}
\end{equation}
For the second equality we used \re{1.3.9}.
On the other hand, in view of the Weingarten equation we obtain
\begin{equation}
\tfrac D{dt}(\nu^\al_i)=\tfrac D{dt}(h^k_i x^\al_k)=\dot h^k_ix^\al_k+h^k_i\dot x^\al_k.
\end{equation}

Multiplying the resulting equation with $\bar g_{\al\bet}x^\bet_j$ we conclude
\begin{equation}
\begin{aligned}
\dot h^k_ig_{kj}-\s(\F-\tilde f)h^k_ih_{kj}&=(\F-\tilde f)_{ij}+\s(\F-\tilde f) \riema \al\bet\ga\de \nu^\al x^\bet_i\nu^\ga x^\de_j
\end{aligned}
\end{equation}
or equivalently \re{1.3.12}.

To derive \re{1.3.13}, we differentiate
\begin{equation}
h_{ij}=h^k_ig_{kj}
\end{equation}
with respect to $t$ and use \re{1.3.6}.  
\ep

We emphasize that equation \re{1.3.12} describes the evolution of the second fundamental form more meaningfully than \re{1.3.13}, since the mixed tensor is independent of the metric.

\bl[Evolution of $(\F-\tilde f)$]\lal{1.3.4}
The term $(\F-\tilde f)$ evolves according to the equation
\begin{equation}\lae{1.3.18}
\begin{aligned}
{(\F-\tilde f)}^\prime-\dot\F F^{ij}(\F-\tilde f)_{ij}=&\msp[3]\sigma \dot \F
F^{ij}h_{ik}h_j^k (\F-\tilde f)
 +\sigma\tilde f_\alpha\n^\alpha (\F-\tilde f)\\
&+\sigma\dot\F F^{ij}\riema \alpha\beta\gamma\delta\n^\alpha x_i^\beta \n^\gamma x_j^\delta (\F-\tilde f),
\end{aligned}
\end{equation}
where
\begin{equation}
(\F-\tilde f)^{\prime}=\frac{d}{dt}(\F-\tilde f)
\end{equation}
and
\begin{equation}
\dot\F=\frac{d}{dr}\F(r).
\end{equation}
\el

\bp
When we differentiate $F$ with respect to $t$ we consider $F$ to depend on the mixed tensor $h^j_i$ and conclude
\begin{equation}
(\F-\tilde f)'=\dot\F F^i_j\dot h^j_i-\tilde f_\al \dot x^\al;
\end{equation}
The equation \re{1.3.18} then follows in view of \re{1.3.5} and \re{1.3.12}.
\ep

\br\lar{1.3.5}
The preceding conclusions, except \rl{1.3.4}, remain valid for flows which do not depend on the curvature, i.e., for flows
\begin{equation}
\begin{aligned}
\dot x&=-\s (-f)\nu=\s f\nu,\\
x(0&)=x_0,
\end{aligned}
\end{equation}
where $f=f(x)$ is defined in an open set $\Om$ containing the initial spacelike hypersurface $M_0$. In the preceding equations we only have to set $\F=0$ and $\tilde f=f$.

The evolution equation for the mean curvature then looks like
\begin{equation}
\dot H=-\D f-\s\{\abs A^2+\bar R_{\al\bet}\nu^\al\nu^\bet\}f,
\end{equation}
where the Laplacian is the Laplace operator on the hypersurface $M(t)$. This is exactly the derivative of the mean curvature operator with respect to normal variations as we shall see in a moment.
\er

But first let us consider the following example.
\bb\lab{1.3.6}
Let $(x^\al)$ be a future directed Gaussian coordinate system in $N$, such that  the metric can be expressed in the form
\begin{equation}
d\bar s^2=e^{2\psi}\{\s(dx^0)^2+\s_{ij}dx^idx^j\}.
\end{equation}
Denote by $M(t)$ the coordinate slices $\{x^0=t\}$, then $M(t)$ can be looked at as the flow hypersurfaces of the flow
\begin{equation}\lae{3.27c}
\dot x=-\s(-e^\psi)\bar\nu,
\end{equation}
where we denote the geometric quantities of the slices by $\bar g_{ij}$, $\bar \nu$, $\bar h_{ij}$, etc. 

Here $x$ is the embedding
\begin{equation}
x=x(t,\xi^i)=(t,x^i).
\end{equation}
Notice that, if $N$ is Riemannian, the coordinate system and the normal are always chosen such that $\nu^0>0$, while, if $N$ is Lorentzian, we always pick the past directed normal.

 Hence the mean curvature of the slices evolves according to
 \begin{equation}\lae{1.3.27}
\dot {\bar H}=-\D e^\psi-\s \{\abs{\bar A}^2+\bar R_{\al\bet}\bar\nu^\al\bar\nu^\bet\}e^\psi.
\end{equation}
\eb

We can now derive the linearization of the mean curvature operator of a spacelike hypersurface, compact or non-compact.

\bn\lan{1.3.7}
Let $M_0\su N$ be a spacelike hypersurface of class $C^4$. We first assume that $M_0$ is compact; then there exists a tubular neighbourhood $\mc U$ and a corresponding \tit{normal} Gaussian coordinate system $(x^\al)$ of class $C^3$ such that $\frac {\pa{}}{\pa x^0}$ is normal to $M_0$. 

Let us consider in 
$\mc U$ of $M_0$ spacelike hypersurfaces $M$ that can be written as graphs
over $M_0$, $M=\graph u$, in the corresponding normal Gaussian coordinate system.
Then the mean curvature of $M$ can be expressed as
\begin{equation}\lae{2.4}
H=\{-\D u+\bar H-\s v^{-2}u^iu^j\bar h_{ij}\}v,
\end{equation}
where $\s=\spd\nu\nu$,  and hence, choosing $u=\e\f$, $\f\in C^2(M_0)$, we deduce
\begin{equation}\lae{2.5}
\begin{aligned}
\frac d{d\e}\fv H {\e=0}&=-\D \f +\dot{\bar H}\f\\[\cma]
&= -\D\f -\s(\abs{\bar A}^2+\bar R_{\al\bet}\nu^\al\nu^\bet)\f,
\end{aligned}
\end{equation}
in view of \re{1.3.27}.

The right-hand side is the \ind{derivative of the mean curvature operator} applied to $\f$. 

If $M_0$ is non-compact,  tubular neighbourhoods exist locally and the relation \re{2.5} will be valid for any $\f\in C^2_c(M_0)$ by using a partition of unity.
\en

The preceding linearization can be immediately generalized to a hypersurface $M_0$ solving the equation
\begin{equation}\lae{3.31c}
\fv F{M_0}=f,
\end{equation}
where $f=f(x)$ is defined in a neighbourhood of $M_0$ and $F=F(h_{ij})$ is curvature operator.

\bl\lal{3.9c}
Let $M_0$ be of class $C^{m,\al}$, $m\ge 2$, $0\le\al\le 1$, satisfy \re{3.31c}. Let $\mc U$ be a (local) tubular neighbourhood of $M_0$, then the linearization of the operator $F-f$ expressed in the normal Gaussian coordinate system $(x^\al)$ corresponding to $\mc U$ and evaluated at $M_0$ has the form
\begin{equation}\lae{3.33c}
-F^{ij}u_{ij}-\s\{F^{ij}h_i^kh_{kj}+F^{ij}\riema \al\bet\ga\de \nu^\al x^\bet_i\nu^\ga x^\de_j+f_\al\nu^\al\}u,
\end{equation}
where $u$ is a function defined in $M_0$, and all geometric quantities are those of $M_0$; the derivatives are covariant derivatives with respect to the induced metric of $M_0$. The operator will be self-adjoint, if $F^{ij}$ is divergence free.
\el

\bp
For simplicity assume that $M_0$ is compact, and let $u\in C^2(M_0)$ be fixed. Then  the hypersurfaces
\begin{equation}
M_\e=\graph(\e u)
\end{equation}
stay in the tubular neighbourhood $\mc U$ for small $\e$, $\abs\e<\e_0$, and their second fundamental forms $(h_{ij})$ can be expressed as
\begin{equation}
v^{-1}h_{ij}=-(\e u)_{ij}+\bar h_{ij},
\end{equation}
where $\bar h_{ij}$ is the second fundamental form of the coordinate slices $\{x^0=\const\}$. 

We are interested in
\begin{equation}
\frac d{d\e}\fv {(F-f)} {\e=0}.
\end{equation}

To differentiate $F$ with respect to $\e$ it is best to consider the mixed form $(h^j_i)$ of the second fundamental form to derive
\begin{equation}
\begin{aligned}
\frac d{d\e}(F-f)=F^i_j\dot h^j_i-\pde f{x^0}u=-F^{ij}u_{ij}+F^i_j\dot {\bar h}^j_i u-\pde f{x^0}u,
\end{aligned}
\end{equation}
where the equation is evaluated at $\e=0$ and $\dot {\bar h}^j_i $ is the derivative of $\bar h^j_i$ with respect to $x^0$.

The result then follows from the evolution equation \re{1.3.12} for the flow \re{3.27c}, i.e., we have to replace $(\F-\tilde f)$ in \re{1.3.12} by $-1$.
\ep

\section{Essential parabolic flow equations}\las{1.4}

From \fre{1.3.12} we deduce with the help of the Ricci identities a parabolic equation
for the second fundamental form
\bl\lal{1.4.1}
The mixed tensor $h_i^j$ satisfies the parabolic equation
\begin{equation}\lae{1.4.1}
\begin{aligned}
&\qq\qq\dot h_i^j-\dot\F F^{kl}h_{i;kl}^j=\\[\cma]
&\hp{=}\;\sigma \dot\F F^{kl}h_{rk}h_l^rh_i^j-\sigma\dot\F F
h_{ri}h^{rj}+\sigma (\F-\tilde f) h_i^kh_k^j\\
&\hp{+}-\tilde f_{\alpha\beta} x_i^\alpha x_k^\beta g^{kj}+\sigma \tilde f_\alpha\n^\alpha h_i^j+\dot\F
F^{kl,rs}h_{kl;i}h_{rs;}^{\hphantom{rs;}j}\\
&\hp{=}+\ddot \F F_i F^j+2\dot \F F^{kl}\riema \alpha\beta\gamma\delta x_m^\alpha x_i ^\beta x_k^\gamma
x_r^\delta h_l^m g^{rj}\\
&\hp{=}-\dot\F F^{kl}\riema \alpha\beta\gamma\delta x_m^\alpha x_k ^\beta x_r^\gamma x_l^\delta
h_i^m g^{rj}-\dot\F F^{kl}\riema \alpha\beta\gamma\delta x_m^\alpha x_k ^\beta x_i^\gamma x_l^\delta h^{mj} \\
&\hp{=}+\sigma\dot\F F^{kl}\riema \alpha\beta\gamma\delta\n^\alpha x_k^\beta\n^\gamma x_l^\delta h_i^j-\sigma\dot\F F
\riema \alpha\beta\gamma\delta\n^\alpha x_i^\beta\n^\gamma x_m^\delta g^{mj}\\
&\hp{=}+\sigma (\F-\tilde f)\riema \alpha\beta\gamma\delta\n^\alpha x_i^\beta\n^\gamma x_m^\delta g^{mj}\\
&\hp{=}+\dot\F F^{kl}\bar R_{\alpha\beta\gamma\delta;\e}\{\n^\alpha x_k^\beta x_l^\gamma x_i^\delta
x_m^\e g^{mj}+\n^\alpha x_i^\beta x_k^\gamma x_m^\delta x_l^\e g^{mj}\}.
\end{aligned}
\end{equation}
\el
\bp
We start with equation \fre{1.3.12} and shall evaluate the term 
\begin{equation}
(\F-\tilde f)^j_i;
\end{equation}
since we are only working with covariant spatial derivatives in the subsequent proof, we may---and shall---consider the covariant form of the tensor
\begin{equation}
(\F-\tilde f)_{ij}.
\end{equation}
First we have
\begin{equation}
\F_i=\dot\F F_i=\dot\F F^{kl}h_{kl;i}
\end{equation}
and
\begin{equation}
\F_{ij}=\dot\F F^{kl}h_{kl;ij}+\Ddot\F F^{kl}h_{kl;i}F^{rs}h_{rs;j}+\dot\F F^{kl,rs}h_{kl,;i}h_{rs;j}.
\end{equation}

Next, we want to replace  $h_{kl;ij}$ by $h_{ij;kl}$. Differentiating the Codazzi equation
\begin{equation}
h_{kl;i}=h_{ik;l}+\bar R_{\al\bet\ga\de}\nu^\al x^\bet_kx^\ga_lx^\de_i,
\end{equation}
where we also used the symmetry of $h_{ik}$, yields
\begin{equation}
\begin{aligned}
&h_{kl;ij}=h_{ik;lj}+\bar R_{\al\bet\ga\de ;\e}\nu^\al x^\bet_kx^\ga_l x^\de_ix^\e_j\\
&+\bar R_{\al\bet\ga\de} \{\nu^\al_j x^\bet_kx^\ga_lx^\de_i+\nu^\al x^\bet_{kj}x^\ga_lx^\de_i+ \nu^\al x^\bet_kx^\ga_{lj}x^\de_i+\nu^\al x^\bet_kx^\ga_lx^\de_{ij}\}.
\end{aligned}
\end{equation}

To replace $h_{kl;ij}$ by $h_{ij;kl}$ we use the Ricci identities
\begin{equation}
h_{ik;lj}=h_{ik;jl}+h_{ak}R^a_{\hp{a}ilj}+h_{ai}R^a_{\hp{a}klj}
\end{equation}
and differentiate once again the Codazzi equation
\begin{equation}
h_{ik;j}=h_{ij;k}+\bar R_{\al\bet\ga\de}\nu^\al x^\bet_ix^\ga_kx^\de_j.
\end{equation}

To replace $\tilde f_{ij}$ we use the chain rule
\begin{equation}
\begin{aligned}
\tilde f_i&=\tilde f_\al x^\al_i,\\
\tilde f_{ij}&=\tilde f_{\al\bet}x^\al_ix^\bet_j+\tilde f_\al x^\al_{ij}.
\end{aligned}
\end{equation}

Then, because of the Gau{\ss} equation, Gaussian formula,  and Weingarten equation,  the symmetry properties of the Riemann curvature tensor and the assumed homogeneity of $F$, i.e.,
\begin{equation}
F=F^{kl}h_{kl},
\end{equation}
we deduce \re{1.4.1} from \fre{1.3.12} after reverting to the mixed representation.
\ep

\br
If we had assumed $F$ to be homogeneous of degree $d_0$ instead of 1, then we
would have to replace the explicit term $F$---occurring twice in the preceding
lemma---by $d_0F$.
\er

If the ambient semi-Riemannian manifold is a space of constant curvature, then the evolution equation of the second fundamental form simplifies considerably, as can be easily verified.

\bl\lal{1.4.3}
Let $N$ be a space of constant curvature $K_N$, then the second fundamental form of the curvature flow \fre{1.3.5} satisfies the  parabolic equation 
\begin{equation}\lae{1.4.12}
\begin{aligned}
\dot h_i^j-\dot\F F^{kl}h_{i;kl}^j&=\s\dot\F F^{kl}h_{rk}h_l^rh_i^j-\s\dot\F F
h_{ri}h^{rj}+\s (\F-\tilde f) h_i^kh_k^j\\
&\hp{=}\;-\tilde f_{\alpha\beta} x_i^\alpha x_k^\beta g^{kj}+\s\tilde f_\alpha\n^\alpha h_i^j+\dot\F
F^{kl,rs}h_{kl;i}h_{rs;}^{\hphantom{rs;}j}\\
&\hp{=}\;+\Ddot\F F_iF^j\\
&\hp{=}\;+K_N\{(\F-\tilde f)\de^j_i+\dot\F F\de^j_i-\dot\F F^{kl}g_{kl}h^j_i\}.
\end{aligned}
\end{equation}
\el

Let us now assume that the open set $\Om\su N$ containing the flow hypersurfaces can be covered by a Gaussian coordinate system $(x^\al)$, i.e., $\Om$ can be topologically viewed as a subset of $I\times \so$, where $\so$ is a compact Riemannian manifold and $I$ an interval. We assume furthermore, that the flow hypersurfaces can be written as graphs over $\so$
\begin{equation}
M(t)=\set{x^0=u(x^i)}{x=(x^i)\in \so};
\end{equation}
we use the symbol $x$ ambiguously by denoting points $p=(x^\al)\in N$ as well as points $p=(x^i)\in \so$ simply by $x$, however, we are careful to avoid confusions.

Suppose that the flow hypersurfaces are given by an embedding $x=x(t,\x)$, where $\x=(\xi^i)$ are local coordinates of a compact manifold $M_0$, which then has to be  homeomorphic to $\so$, then \begin{equation}
\begin{aligned}
x^0&=u(t,\xi)=u(t,x(t,\xi)),\\
x^i&=x^i(t,\xi).
\end{aligned}
\end{equation}

The induced metric can be expressed as
\begin{equation}
\begin{aligned}
g_{ij}=\spd {x_i}{x_j}=\s u_iu_j+\s_{kl}x^k_ix^l_j,
\end{aligned}
\end{equation}
where
\begin{equation}
u_i=u_kx^k_i,
\end{equation}
i.e.,
\begin{equation}
g_{ij}=\{\s u_ku_l+\s_{kl}\}x^k_ix^l_j,
\end{equation}
hence the (time dependent)  \tit{\ind{Jacobian}} $(x^k_i)$ is invertible, and the $(\xi^i)$ can also be viewed as coordinates for $\so$.

Looking at the component $\al=0$ of the flow equation \fre{1.3.5} we obtain a scalar flow equation
\begin{equation}\lae{1.4.18}
\dot u=-e^{-\psi}v^{-1}(\F-\tilde f),
\end{equation}
which is the same in the Lorentzian as well as in the Riemannian case, where
\begin{equation}
v^2=1-\s \s^{ij}u_iu_j,
\end{equation}
and where 
\begin{equation}
\abs{Du}^2=\s^{ij}u_iu_j
\end{equation}
is of course a scalar, i.e., we obtain the same expression regardless, if we use the coordinates $x^i$ or $\xi^i$.

The time derivative in \re{1.4.18} is a total time derivative, if we consider $u$ to depend on $u=u(t,x(t,\xi))$. For the partial time derivative we obtain
\begin{equation}\lae{1.4.21}
\begin{aligned}
\frac {\pa u}{\pa t}&=\dot u-u_k\dot x^k_i\\
&=-e^{-\psi}v(\F-\tilde f),
\end{aligned}
\end{equation}
in view of \fre{1.3.5} and our choice of normal $\n=(\nu^\al)$
\begin{equation}
(\nu^\al)=\s e^{-\psi}v^{-1}(1,-\s u^i),
\end{equation}
where $u^i=\s^{ij}u_j$.

Controlling the $C^1$-norm of the graphs $M(t)$ is tantamount to controlling $v$, if $N$ is Riemannian, and $\tilde v=v^{-1}$, if $N$ is Lorentzian. The evolution equations satisfied by these quantities are also very important, since they are used for the a priori estimates of the second fundamental form.

Let us start with the Lorentzian case.

\bl[Evolution of $\tilde v$]\lal{1.4.4}
Consider the flow \re{1.3.5} in a Lorentzian space $N$ such that the spacelike flow hypersurfaces can be written  as graphs over $\so$. Then, $\tilde v$ satisfies the evolution equation
\begin{equation}\lae{1.4.23}
\begin{aligned}
\dot{\tilde v}-\dot\F F^{ij}\tilde v_{ij}=&-\dot\F F^{ij}h_{ik}h_j^k\tilde v
+[(\F-\tilde f)-\dot\F F]\h_{\alpha\beta}\n^\alpha\n^\beta\\
&-2\dot\F F^{ij}h_j^k x_i^\alpha x_k^\beta \h_{\alpha\beta}-\dot\F F^{ij}\h_{\alpha\beta\gamma}x_i^\beta
x_j^\gamma\n^\alpha\\
&-\dot\F F^{ij}\riema \alpha\beta\gamma\delta\n^\alpha x_i^\beta x_k^\gamma x_j^\delta\h_\e x_l^\e g^{kl}\\
&-\tilde f_\beta x_i^\beta x_k^\alpha \h_\alpha g^{ik},
\end{aligned}
\end{equation}
where $\h$ is the covariant vector field $(\h_\alpha)=e^{\psi}(-1,0,\dotsc,0)$.
\el

\bp
We have $\tilde v=\spd \h\n$. Let $(\x^i)$ be local coordinates for $M(t)$.
Differentiating $\tilde v$ covariantly we deduce
\begin{equation}\lae{1.4.24}
\tilde v_i=\h_{\alpha\beta}x_i^\beta\n^\alpha+\h_\alpha\n_i^\alpha,
\end{equation}
\begin{equation}\lae{1.4.25}
\begin{aligned}
\tilde v_{ij}= &\msp[5]\h_{\alpha\beta\gamma}x_i^\beta x_j^\gamma\n^\alpha+\h_{\alpha\beta}x_{ij}^\beta\n^\alpha\\
&+\h_{\alpha\beta}x_i^\beta\n_j^\alpha+\h_{\alpha\beta}x_j^\beta\n_i^\alpha+\h_\alpha\n_{ij}^\alpha
\end{aligned}
\end{equation}

The time derivative of $\tilde v$ can be expressed as
\begin{equation}\lae{1.4.26}
\begin{aligned}
\dot{\tilde v}&=\h_{\alpha\beta}\msp\dot x^\beta\n^\alpha+\h_\alpha\dot\n^\alpha\\
&=\h_{\alpha\beta}\n^\alpha\n^\beta(\F-\tilde f)+(\F-\tilde f)^k x_k^\alpha\h_\alpha\\
&=\h_{\alpha\beta}\n^\alpha\n^\beta(\F-\tilde f)+\dot\F F^k x_k^\alpha\h_\alpha-{\tilde f}_\beta x_i^\beta
x_k^\alpha g^{ik}\h_\alpha,
\end{aligned}
\end{equation}
where we have used \fre{1.3.9}.

Substituting \re{1.4.25} and \re{1.4.26} in \re{1.4.23}, and simplifying the resulting
equation with the help of the Weingarten and Codazzi equations, we arrive at the
desired conclusion.
\ep

In the Riemannian case we consider a normal Gaussian coordinate system $(x^\al)$, for otherwise we won't obtain a priori estimates for $v$, at least not without additional strong assumptions. We also refer to $x^0=r$ as the \tit{\ind{radial distance function}}.

\bl[Evolution of $v$]\lal{1.4.5}
Consider the flow \re{1.3.5} in a normal Gaussian coordinate system where the $M(t)$ can be written as graphs of a function $u(t)$ over some compact Riemannian manifold $\so$. Then the quantity
\begin{equation}\lae{1.4.28}
v=\sqrt{1+\abs {Du}^2}=(r_\al\nu^\al)^{-1}
\end{equation}
satisfies the evolution equation
\begin{equation}
\begin{aligned}
&\dot v -\dot\F F^{ij}v_{ij}= -\dot\F F^{ij}h_{ik}h^k_jv-2v^{-1}\dot\F F^{ij}v_iv_j\\
&\hp{=}\q +r_{\al\bet}\nu^\al\nu^\bet[(\F-\tilde f)-\dot\F F]v^2+2\dot\F F^{ij}h^k_ir_{\al\bet}x^\al_k x^\bet_j v^2\\
&\hp{=}\q +\dot\F F^{ij}\bar R_{\al\bet\ga\de}\nu^\al x^\bet_i x^\ga_jx^\de_k r_\e x^\e_m g^{mk}v^2\\
&\hp{=}\q +\dot \F F^{ij}r_{\al\bet\ga}\nu^\al x^\bet_i x^\ga _j v^2 +\tilde f_\al x^\al_mg^{mk}r_\bet x^\bet_k v^2.
\end{aligned}
\end{equation}
\el

\bp
Similar to the proof of the previous lemma. 
\ep


The previous problems can be generalized to the case when the right-hand side $f$ is not only defined in $N$ or in $\bar\Om$ but in the \tit{\ind{tangent bundle}} $T(N)$ \resp $T(\bar\Om)$. Notice that the tangent bundle is a manifold of dimension $2(n+1)$, i.e., in a \tit{\ind{local trivialization}} of $T(N)$ $f$ can be expressed in the form
\begin{equation}
f=f(x,\nu)
\end{equation}
with $x\in N$ and $\nu\in T_x(N)$, \cf \cite[Note 12.2.14]{cg:aII}. Thus, the case $f=f(x)$ is included in this general set up. The symbol $\nu$ indicates that in an equation
\begin{equation}\lae{1.4.30}
\fv FM=f(x,\nu)
\end{equation}
we want $f$ to be evaluated at $(x,\nu)$, where $x\in M$ and $\nu$ is the normal of $M$ in $x$.

The \tit{\ind{Minkowski problem}} or \tit{\ind{Minkowski type problems}} are also covered by the present setting, though the Minkowski problem has the additional property that the problem is transformed via the \tit{\ind{Gau{\ss} map}} to a different semi-Riemannian manifold as a \tit{\ind{dual problem}} and solved there. Minkowski type problems have been treated in \cite{cheng-yau:minkowski}, \cite{guan:annals}, \cite{cg:minkowski} and \cite{cg:minkowski2}.
\br
 The equation \re{1.4.30} will be solved by the same methods as in the special case when $f=f(x)$, i.e., we consider the same curvature flow, the evolution equation \fre{1.3.5}, as before.
 
 The resulting evolution equations are identical with the natural exception, that, when $f$ or $\tilde f$ has to be differentiated, the additional argument has to be considered, e.g.,
 \begin{equation}\lae{1.4.31}
\tilde f_i=\tilde f_\al x^\al_i+\tilde f_{\nu^\bet}\nu^\bet_i=\tilde f_\al x^\al_i+\tilde f_{\nu^\bet}x^\bet_k h^k_i
\end{equation}
and
\begin{equation}\lae{1.4.32}
\dot {\tilde f}=\tilde f_\al \dot x^\al+\tilde f_{\nu^\bet}\dot\nu^\bet=-\s(\F-\tilde f)\tilde f_\al\nu^\al+\tilde f_{\nu^\bet}g^{ij}(\F-\tilde f)_i x^\bet_j.
\end{equation}

The most important evolution equations are explicitly stated below. 
\er

Let us first state the evolution equation for $(\F-\tilde f)$.

\bl[Evolution of $(\F-\tilde f)$]\lal{1.4.7}
The term $(\F-\tilde f)$ evolves according to the equation
\begin{equation}\lae{1.4.3.10}
\begin{aligned}
{(\F-\tilde f)}^\prime-\dot\F F^{ij}(\F-\tilde f)_{ij}&=
\s \dot \F
F^{ij}h_{ik}h_j^k (\F-\tilde f)\\
&\hp{=}\;+\s\tilde f_\al\n^\al (\F-\tilde f)
- \tilde f_{\n^\al}x^\al_i(\F- \tilde
f)_jg^{ij}\\
&\hp{=}\;+\s\dot\F F^{ij}\riema \al\bet\ga\de\n^\al x_i^\bet \n^\ga x_j^\de
(\F-\tilde f),
\end{aligned}
\end{equation}where
\begin{equation}
(\F-\tilde f)^{\prime}=\frac{d}{dt}(\F-\tilde f)
\end{equation}
and
\begin{equation}
\dot\F=\frac{d}{dr}\F(r).
\end{equation}
\el

Here is the evolution equation for the second fundamental form.
\bl\lal{1.4.8}
The mixed tensor $h_i^j$ satisfies the parabolic equation
\begin{equation}\lae{1.4.3.13}
\begin{aligned}
\dot h_i^j&-\dot\F F^{kl}h_{i;kl}^j\\[\cma]&=\s \dot\F
F^{kl}h_{rk}h_l^rh_i^j-\s\dot\F F h_{ri}h^{rj}+\s (\F-\tilde f)
h_i^kh_k^j\\[\cma] 
&\hp{+}-\tilde f_{\al\bet} x_i^\al x_k^\bet g^{kj}+\s \tilde f_\al\n^\al
h_i^j-\tilde f_{\al\n^\bet}(x^\al_i x^\bet_kh^{kj}+x^\al_l x^\bet_k h^{k}_i\, g^{lj})\\
&\hp{=}
-\tilde f_{\n^\al\n^\bet}x^\al_lx^\bet_kh^k_ih^{lj}-\tilde f_{\n^\bet} x^\bet_k
h^k_{i;l}\,g^{lj}  +\s\tilde f_{\n^\al}\n^\al h^k_i h^j_k\\
&\hp{=}+\dot\F
F^{kl,rs}h_{kl;i}h_{rs;}^{\hphantom{rs;}j}+2\dot \F F^{kl}\riema \al\bet\ga\de x_m^\al x_i ^\bet x_k^\ga
x_r^\de h_l^m g^{rj}\\
&\hp{=}-\dot\F F^{kl}\riema \al\bet\ga\de x_m^\al x_k ^\bet x_r^\ga x_l^\de
h_i^m g^{rj}-\dot\F F^{kl}\riema \al\bet\ga\de x_m^\al x_k ^\bet x_i^\ga x_l^\de h^{mj} \\
&\hp{=}+\s\dot\F F^{kl}\riema \al\bet\ga\de\n^\al x_k^\bet\n^\ga x_l^\de h_i^j-\s\dot\F F
\riema \al\bet\ga\de\n^\al x_i^\bet\n^\ga x_m^\de g^{mj}\\
&\hp{=}+\s (\F-\tilde f)\riema \al\bet\ga\de\n^\al x_i^\bet\n^\ga x_m^\de g^{mj}+\ddot \F
F_i F^j\\ &\hp{=}+\dot\F F^{kl}\bar R_{\al\bet\ga\de;\e}\{\n^\al x_k^\bet x_l^\ga x_i^\de
x_m^\e g^{mj}+\n^\al x_i^\bet x_k^\ga x_m^\de x_l^\e g^{mj}\}.
\end{aligned}
\end{equation}
\el

The proof is identical to that of  \rl{1.4.1};
we only have to keep in mind that $f$ now also depends on the normal.

If we had assumed $F$ to be homogeneous of degree $d_0$ instead of $1$, then, we
would have to replace the explicit term $F$---occurring twice in the preceding
lemma---by $d_0F$.

\bl[Evolution of $\tilde v$]\lal{1.4.9}
Consider the flow \re{1.3.5} in a Lorentzian space $N$ such that the spacelike flow hypersurfaces can be written  as graphs over $\so$. Then, $\tilde v$ satisfies the evolution equation
\begin{equation}\lae{1.4.4.3}
\begin{aligned}
\dot{\tilde v}-\dot\F F^{ij}\tilde v_{ij}=&-\dot\F F^{ij}h_{ik}h_j^k\tilde v
+[(\F-\tilde f)-\dot\F F]\h_{\al\bet}\n^\al\n^\bet\\
&-2\dot\F F^{ij}h_j^k x_i^\al x_k^\bet \h_{\al\bet}-\dot\F F^{ij}\h_{\al\bet\ga}x_i^\bet
x_j^\ga\n^\al\\
&-\dot\F F^{ij}\riema \al\bet\ga\de\n^\al x_i^\bet x_k^\ga x_j^\de\h_\e x_l^\e g^{kl}\\
&-\tilde f_\bet x_i^\bet x_k^\al \h_\al g^{ik}- \tilde f_{\n^\bet}x^\bet_k h^{ik}x^\al_i\h_\al,
\end{aligned}
\end{equation}
where $\h$ is the covariant vector field $(\h_\al)=e^{\psi}(-1,0,\dotsc,0)$.
\el
The proof is identical to the proof of \rl{1.4.4}.

In the Riemannian case we have:

\bl[Evolution of $v$]\lal{1.4.10}
Consider the flow \re{1.3.5} in a normal Gaussian coordinate system $(x^\al)$, where the $M(t)$ can be written as graphs of a function $u(t)$ over some compact Riemannian manifold $\so$. Then the quantity
\begin{equation}
v=\sqrt{1+\abs {Du}^2}=(r_\al\nu^\al)^{-1}
\end{equation}
satisfies the evolution equation
\begin{equation}\lae{1.4.4.10}
\begin{aligned}
\dot{ v}-\dot\F F^{ij} v_{ij}=&-\dot\F F^{ij}h_{ik}h_j^k v-2 v^{-1} \dot\F
F^{ij}v_iv_j \\
&+[(\F- f)-\dot\F F]r_{\al\bet}\n^\al\n^\bet v^2\\
&+2\dot\F F^{ij}h_j^k x_i^\al x_k^\bet r_{\al\bet} v^2+\dot\F F^{ij}r_{\al\bet\ga}x_i^\bet
x_j^\ga\n^\al v^2\\
&+\dot\F F^{ij}\riema \al\bet\ga\de\n^\al x_i^\bet x_k^\ga x_j^\de r_\e x_l^\e g^{kl} v^2\\
&+\tilde f_\bet x_i^\bet x_k^\al r_\al g^{ik}v^2+\tilde f_{\n^\bet}x^\bet_k h^{ik}x^\al_i
r_\al v^2,
\end{aligned}
\end{equation}
where $r=x^0$ and $(r_\al)=(1,0,\dotsc,0)$.
\el

\section{Stability of the limit hypersurfaces}\las{5c}

\bd
Let $N$ be semi-Riemannian,  $F$ a curvature operator, and $M\su N$ a compact, spacelike hypersurface, such that $M$ is admissible and $F^{ij}$, evaluated at $(h_{ij},g_{ij})$, the second fundamental form and metric of $M$, is divergence free, then $M$ is said to be a \tit{stable} solution of the equation
\begin{equation}
\fv FM=f,
\end{equation}
where $f=f(x)$ is defined in a neighbourhood of $M$, if the first eigenvalue $\lam_1$ of the linearization, which is the operator in \fre{3.33c}, is non-negative, or equivalently, if the quadratic form
\begin{equation}
\int_M F^{ij}u_iu_j-\s\int_M\{F^{ij}h_i^kh_{kj}+F^{ij}\riema \al\bet\ga\de \nu^\al x^\bet_i\nu^\ga x^\de_j+f_\al\nu^\al\}u^2
\end{equation}
is non-negative for all $u\in C^2(M)$.
\ed

It is well-known that the corresponding eigenspace is then onedimensional and spanned by a strictly positive eigenfunction $\h$
\begin{equation}\lae{5.3c}
-F^{ij}\h_{ij}-\s\{F^{ij}h_i^kh_{kj}+F^{ij}\riema \al\bet\ga\de \nu^\al x^\bet_i\nu^\ga x^\de_j+f_\al\nu^\al\}\h=\lam_1\h.
\end{equation}
Notice that $F^{ij}$ is supposed to be divergence free, which will be the case, if $F=H_k$, $1\le k\le n$, and the ambient space has constant curvature, as we shall prove at the end of this section. If $k=1$, then $F^{ij}=g^{ij}$ and $N$ can be arbitrary, while in case $k=2$, we have
\begin{equation}
F^{ij}=Hg^{ij}-h^{ij},
\end{equation}
hence $N$ Einstein will suffice.

To simplify the formulation of the assumptions let us define:
\bd
A curvature function $F$ is said to be of class $(D)$, if for every admissible hypersurface $M$ the tensor $F^{ij}$, evaluated at $M$, is divergence free.
\ed

We shall prove in this section that the limit hypersurface of a converging curvature flow will be a stable stationary solution, if the initial flow velocity has a weak sign.

\bt\lat{5.3c}
Suppose that the curvature flow \fre{1.3.5} exists for all time, and that the leaves $M(t)$ converge in $C^4$ to a hypersurface $M$, where the curvature function $F$ is supposed to be of class $(D)$. Then $M$ is a stable solution of the equation
\begin{equation}\lae{5.5c}
\fv FM=f
\end{equation}
provided the velocity of the flow has a weak sign
\begin{equation}\lae{5.6c}
\F-\tilde f\ge0\q\vee\q \F-\tilde f\le0
\end{equation}
at $t=0$ and $M(0)$ is not already a solution of \re{5.5c}.
\et

\bp
Convergence of a subsequence of the $M(t)$ would actually suffice for the proof, however, the assumption \re{5.6c} immediately implies that the flow converges, if a subsequence converges and a priori estimates in $C^{4,\al}$ are valid. 

The starting point is the evolution equation \fre{1.3.18} from which we deduce in view of the parabolic maximum principle that $\F-\tilde f$ has a weak sign during the evolution, \cf \cite[Proposition 2.7.1]{cg:cp}, i.e., if we assume without loss of generality that at $t=0$
\begin{equation}
\F-\tilde f\ge0,
\end{equation}
then this inequality will be valid for all $t$. Moreover, there holds
\begin{equation}\lae{5.8a}
\int_{M(t)}(\F-\tilde f)>0\qq\A\,0\le t<\un
\end{equation}
if this relation is valid for $t=0$, as we shall prove in the lemma below. 

On the other hand, the assumption
\begin{equation}
\int_{M(0)}(\F-\tilde f)>0
\end{equation}
 is a natural assumption, for otherwise the initial hypersurface would already be a stationary solution which of course may not be stable.

Notice also that apart from the factor $\dot\F$ the equation \re{1.3.18} looks like the parabolic version of the linearization of $(F-f)$. If the technical function $\F=\F(r)$ is not the trivial one $\F(r)=r$, then we always assume that $f>0$ and that this is also valid for the limit hypersurface $M$. Only in case $\F(r)=r$ and $F=H$, we allow $f$ to be arbitrary.

Thus, our assumptions imply  that in any case
\begin{equation}
\dot\F>\e_0>0\qq\A\, t\in\R[]_+.
\end{equation}

Furthermore, we derive from \re{1.3.18} that not only the elliptic part converges to $0$ but also 
\begin{equation}
(\F-\tilde f)'=\dot\F \dot F +\s\dot\F(f) f_\al\nu^\al(\F-\tilde f),
\end{equation}
i.e.,
\begin{equation}\lae{5.10c}
\lim \dot F=0.
\end{equation}

Suppose now that $M$ is not stable, then the first eigenvalue $\lam_1$ is negative and there exists a strictly positive eigenfunction $\h$ solving the equation \re{5.3c} evaluated at $M$. Let $\mc U$ be a tubular neighbourhood of $M$ with a corresponding future directed normal Gaussian coordinate system $(x^\al)$ and extend $\h$ to $\mc U$ by setting
\begin{equation}
\h(x^0,x)=\h(x),
\end{equation}
where, by a slight abuse of notation, we also denote $(x^i)$ by $x$. Thus there holds
\begin{equation}
\h_\al\nu^\al=0
\end{equation}
in $M$, and choosing $\mc U$ sufficiently small, we may assume
\begin{equation}\lae{5.13c}
\abs{\h_\al\nu^\al}<\h
\end{equation}
for all hypersurfaces $M(t)\su\mc U$.

Now consider the term
\begin{equation}
\int_{M(t)}\dot\F^{-1} (\F-\tilde f)\h
\end{equation}
for large $t$, which converges to $0$. Since it is positive, in view of \re{5.8a}, there must exist a sequence of $t$, not explicitly labelled, tending to infinity such that
\begin{equation}\lae{5.15c}
\begin{aligned}
0&\ge \frac d{dt}\int_{M(t)}\dot\F^{-1} (\F-\tilde f)\h\\
&=\int_{M(t)}-\dot\F^{-2}\Ddot\F\dot F(\F-\tilde f)\h
+\int_{M(t)}\dot\F^{-1}(\F-\tilde f)'\h\\
&\hp{=}-\s\int_{M(t)}\dot\F^{-1}\h_\al\nu^\al(\F-\tilde f)^2-\s\int_{M(t)}\dot\F^{-1}(\F-\tilde f)^2 H\h,
\end{aligned}
\end{equation}
where we used the relation \fre{1.3.6} to derive the last integral.

The rest of the proof is straight-forward. Multiply the equation \re{1.3.18} by $\dot\F^{-1}\h$ and integrate over $M(t)$ for those values of $t$ satisfying the preceding inequality to deduce
\begin{equation}
\begin{aligned}
&\qq -\int_{M(t)}\dot\F^{-1}(\F-\tilde f)'=\int_{M(t)}-F^{ij}\h_{ij}(\F-\tilde f)\\
&-\s\int_{M(t)}\{F^{ij}h_i^kh_{kj}+F^{ij}\riema \al\bet\ga\de \nu^\al x^\bet_i\nu^\ga x^\de_j+\dot\F^{-1}\dot\F(f) f_\al\nu^\al\}\h(\F-\tilde f),
\end{aligned}
\end{equation}
and conclude further that the right-hand side can be estimated from above by 
\begin{equation}
\frac{\lam_1}2\int_{M(t)}\h(\F-\tilde f)
\end{equation}
for large $t$, while the left-hand side can be estimated from below by
\begin{equation}
-\e(t)\int_{M(t)}(\F-\tilde f)\h
\end{equation}
such that 
\begin{equation}
\e(t)>0\q\wed\q\lim\e(t)=0
\end{equation}
in view of \re{5.15c}, where we used \re{5.10c}, \re{5.13c}  as well as
\begin{equation}
\lim(\F-\tilde f)=0;
\end{equation}
a contradiction because of \re{5.8a}.
\ep

\bl\lal{5.4c}
Let $M(t)$ be a solution of the curvature flow \fre{1.3.5} defined on a maximal time interval $[0,T^*)$, $0<T^*\le \un$, and suppose that $\F-\tilde f$ has a weak sign at $t=0$, e.g.,
\begin{equation}
(\F-\tilde f)\ge0
\end{equation}
and suppose furthermore that
\begin{equation}
\int_{M(0)}(\F-\tilde f)>0,
\end{equation}
then
\begin{equation}
\int_{M(t)}(\F-\tilde f)>0\qq\A\, 0\le t<T^*.
\end{equation}
\el

\bp
Let $M_0$ be an abstract compact Riemannian manifold that is being isometrically embedded in $N$ with image $M(0)$. Let $(\xi^i)$ be a generic coordinate system for $M_0$ and abbreviate $(\F-\tilde f)$ by $u$. The evolution equation \re{1.3.18} can then be looked at as a linear parabolic equation for $u=u(t,\xi)$ on $M_0$ with time dependent coefficients and time dependent Riemannian metric $g_{ij}(t,\xi)$.

By assumption $u$ doesn't vanish identically at $t=0$, i.e., there exists a ball $B_\rho=B_\rho(\xi_0)$ such that 
\begin{equation}
u(0,\xi)>0\qq\A\,\xi\in \bar B_\rho(\xi_0).
\end{equation}
Let $\mc C$ be the cylinder
\begin{equation}
\mc C=[0,T^*)\times \bar B_\rho
\end{equation}
and assume that there exists a first $t_0>0$ such that
\begin{equation}
\inf_{\bar B_\rho}u(t_0,\cdot)=0=u(t_0, \xi_1).
\end{equation}
We shall show that this is not possible: If $\xi_1\in B_\rho$, then this contradicts the strong parabolic maximum principle, \cf \cite[Lemma 2.7.1]{cg:cp}, and if $\xi_1\in \pa B_\rho$, then we deduce from
\cite[Lemma 2.7.4]{cg:cp} (a parabolic version of the Hopf Lemma)
\begin{equation}
\pde u{\nu}(t_0,\xi_1)<0,
\end{equation}
where $\nu$ is the exterior normal of the ball $B_\rho$ in $\xi_1$, contradicting the fact that the gradient of $u(t_0,\cdot)$ vanishes in $\xi_1$ because it is a minimum point; notice that we already know $u\ge 0$ in $[0,T^*)\times M_0$.
\ep

For some curvature operators one can prove a priori estimates for the second fundamental form only for stationary solutions and not for the leaves of a corresponding curvature flow. In order to use a curvature flow to obtain a stationary solution one uses \cqp{$\e$-regularization}, i.e., instead of the curvature function $F$ one considers
\begin{equation}\lae{5.31.1}
\tilde F(h_{ij})=F(h_{ij}+\e H g_{ij})
\end{equation}
for $\e>0$, and starts a curvature flow with $\tilde F$ and fixed $\e>0$.

A priori estimates for the regularized flow are usually fairly easily derived, since
\begin{equation}
\tilde F^{ij}=F^{ij}+\e F^{kl}g_{kl}g^{ij},
\end{equation}
but of course the estimates depend on $\e$. Having uniform estimates one can deduce that the flow---or at least a subsequence---converges to a limit hypersurfaces $M_\e$  satisfying
\begin{equation}
\fv{\tilde F}{M_\e}=f.
\end{equation}
Then, if uniform $C^{4,\al}$-estimates for the $M_\e$ can be derived, a subsequence will converge to a solution $M$ of
\begin{equation}
\fv FM=f,
\end{equation}
\cf \cite{cg:scalar}, where this method has been used to find hypersurfaces of prescribed scalar curvature in Lorentzian manifolds, see also \frt{6.6.1}. 

We shall now show that the solutions $M$ obtained by this approach are all stable, if $F$ is of class $(D)$ and the initial velocities of the regularized flows have a weak sign. Notice that the curvature functions $\tilde F$ are in general not of class $(D)$.

\bt
Let $F$ be of class $(D)$, then any solution $M$ of
\begin{equation}
\fv FM=f
\end{equation}
obtained by a regularized curvature flow as described above is stable, provided the initial velocity of the regularized flow has a weak sign, i.e., it satisfies
\begin{equation}
\F-\tilde f\ge0\q\vee\q \F-\tilde f\le 0
\end{equation}
at $t=0$ and the flow hypersurfaces converge to the stationary solution in $C^4$.
\et

\bp
Let $M_\e$ be the limit hypersurfaces of the regularized flow for $\e>0$, and assume that the $M_\e$ satisfy uniform $C^{4,\al}$-estimates such that a subsequence, not relabelled, converges in $C^4$ to a compact spacelike hypersurface $M$ solving the equation
\begin{equation}
\fv FM=f.
\end{equation}

Assume that $M$ is not stable so that the first eigenvalue of the linearization is negative and there exists a strictly positive eigenfunction $\h$ satisfying \re{5.3c}. Extend $\h$ in a small tubular neighbourhood $\mc U$ of $M$ such that \re{5.13c} is valid for all $M_\e$, if $\e$ is small, $\e<\e_0$.

For those $\e$ we then deduce
\begin{equation}
\begin{aligned}
-\tilde F^{ij}\h_{ij}&-\tilde F^{ij}_{;ij}\h-2\tilde F^{ij}_{;j}\h_i\\
&-\s\{\tilde F^{ij}h_i^kh_{kj}+\tilde F^{ij}\riema \al\bet\ga\de \nu^\al x^\bet_i\nu^\ga x^\de_j+f_\al\nu^\al\}\h<\frac{\lam_1}2\h,
\end{aligned}
\end{equation}
where the inequality is evaluated at $M_\e$ and where we used the convergence in $C^4$.

Now, fix $\e$, $\e<\e_0$, then the preceding inequality is also valid for the flow hypersurfaces $M(t)$ converging to $M_\e$, if $t$ is large, and the same arguments as at the end of the proof of \rt{5.3c} lead to a contradiction. Hence, $M$ has to be a stable solution.
\ep

Knowing that a solution is stable often allows to deduce further geometric properties of the underlying hypersurface like that it is either strictly stable or totally geodesic especially if the curvature function is the mean curvature, \cf e.g., \cite{yau:stable}, where the stability property has been extensively used to deduce geometric properties.

We want to prove that a neighbourhood of stable solutions can be foliated by a family of hypersurfaces satisfying the equation modulo a constant.

\bt
Let $M\su N$ be compact, spacelike, orientable and  a stable solution of
\begin{equation}
\fv GM\equiv \fv {(F-f)}M=0,
\end{equation}
where $F$ is of class $(D)$ and $M$ as well as $F$, $f$ are of class $C^{m,\al}$, $2\le m\le\un$, $0<\al<1$, then a neighbourhood of $M$ can be foliated by a family
\begin{equation}
\Lam=\set{M_\e}{\abs \e<\e_0}
\end{equation}
of spacelike $C^{m,\al}$-hypersurfaces satisfying
\begin{equation}
\fv G{M_\e}=\tau(\e),
\end{equation}
where $\tau$ is a real function of class $C^{m,\al}$. The $M_\e$ can be written as graphs over $M$ in a tubular neighbourhood of $M$
\begin{equation}
M_\e=\set{(u(\e,x),x)}{x\in M}
\end{equation}
such that $u$ is of class $C^{m,\al}$ in both variables and there holds
\begin{equation}
\dot u>0.
\end{equation}
\et

\bp
(i) Let us assume that $M$ is strictly stable. Consider a tubular neighbourhood of $M$ with corresponding normal Gaussian coordinates $(x^\al)$ such that $M=\{x^0=0\}$. The nonlinear operator $G$ can then be viewed as an elliptic operator
\begin{equation}
G:B_\rho(0)\su C^{m,\al}(M)\ra C^{m-2,\al}(M)
\end{equation}
where $\rho$ is so small that all corresponding graphs are admissible.

In a smaller ball $DG$ is a topological isomorphism, since $M$ is strictly stable, and hence $G$ is a diffeomorphism in a neighbourhood of the origin, and there exist smooth unique solutions
\begin{equation}
M_\e=\set{u(\e,x)}{x\in M}\qq\abs\e<\e_0
\end{equation}
of the equations
\begin{equation}
\fv G{M_\e}=\e
\end{equation}
such that $u\in C^{m,\al}((-\e_0,\e_0)\times M)$.

Differentiating with respect to $\e$ yields
\begin{equation}\lae{5.46c}
DG\dot u=1.
\end{equation}

Let us consider this equation for $\e=0$, i.e., on $M$, and define
\begin{equation}
\h=\min(\dot u,0).
\end{equation}
Then we deduce
\begin{equation}
0\le \int_M\spd{DG\h}{\h}=\int_M\h\le0,
\end{equation}
and hence there holds
\begin{equation}
\dot u\ge 0,
\end{equation}
because of the strict stability of $M$.

Applying then the maximum principle to \re{5.46c}, we deduce further
\begin{equation}
\inf_M\dot u>0,
\end{equation}
hence the hypersurfaces form a foliation if $\e_0$ is chosen small enough such that
\begin{equation}
\inf_{M_\e}\dot u(\e,\cdot)>0\qq\A\, \abs\e<\e_0.
\end{equation}

\cvm
(ii) Assume now that $M$ is not strictly stable. After introducing coordinates corresponding to a tubular neighbourhood $\mc U$ of $M$ as in part (i) any function $u\in C^{m,\al}(M)$ with $\abs u_{m,\al}$ small enough defines an admissible hypersurface
\begin{equation}
M(u)=\graph u\su\mc U
\end{equation}
such that $\fv G{M(u)}$ can be expressed as
\begin{equation}
\fv G{M(u)}=G(u).
\end{equation}

Let
\begin{equation}
A=DG(0),
\end{equation}
then $A$ is self-adjoint, monotone
\begin{equation}
\spd{Au}u\ge 0\qq\A\, u\in H^{1,2}(M)
\end{equation}
and the smallest eigenvalue of $A$ is equal to zero, the corresponding eigenspace spanned by a strictly positive eigenfunction $\h$.

Similarly as in \cite[p. 621]{bartnik:1988} we consider the operator
\begin{equation}
\Psi(u,\tau)=(G(u)-\tau,\f(u))
\end{equation}
defined in $B_\rho(0)\times\R[]$, $B_\rho(0)\su C^{m,\al}(M)$ for small $\rho>0$, where $\f$ is a linear functional
\begin{equation}
\f(u)=\int_M\h u
\end{equation}
$\Psi$ is of class $C^{m,\al}$ and maps
\begin{equation}
\Psi:B_\rho(0)\times \R[]\ra C^{m-2,\al}(M)\times \R[],
\end{equation}
such that
\begin{equation}
D\Psi=
\begin{pmatrix}
DG&-1\\[\cma]
\f&0
\end{pmatrix}
\end{equation}
evaluated at $(0,0)$ is bijective  as one easily checks. Indeed let $(u,\e)$ satisfy
\begin{equation}
D\Psi(u,\e)=(0,0),
\end{equation}
then
\begin{equation}
Au=DGu=\e\q\wed\q\int_M\h u=0,
\end{equation}
hence
\begin{equation}
\e\int_M\h=\spd{Au}\h=\spd u{A\h}=0
\end{equation}
and we conclude $\e=0$ as well as $u=0$.

To prove the surjectivity, let $(w,\de)\in C^{m-2,\al}(M)\times \R[]$ be arbitrary. Choosing
\begin{equation}
\e=-\frac{\int_Mw\h}{\int_M\h}
\end{equation}
we deduce
\begin{equation}
\int_M(\e+w)\h=0,
\end{equation}
hence there exists $\bar u\in C^{m,\al}(M)$ solving
\begin{equation}
A\bar u=\e+w
\end{equation}
and 
\begin{equation}
u=\bar u+\lam\h
\end{equation}
with
\begin{equation}
\lam=\de-\int_M\h\bar u
\end{equation}
then satisfies
\begin{equation}
\int_M\h u=\e,
\end{equation}
i.e.,
\begin{equation}
D\Psi(u,\e)=(w,\de).
\end{equation}

Applying the inverse function theorem we conclude that  there exists $\e_0>0$ and functions $(u(\e,x),\tau(\e))$ of class $C^{m,\al}$ in both variables such that
\begin{equation}\lae{5.70c}
G(u(\e))=\tau(\e)\q\wed\q\int_M\h u(\e)=\e\q\A\, \abs\e<\e_0;
\end{equation}
 $\tau(\e)$ is  constant for fixed $\e$.
 
 The hypersurfaces
 \begin{equation}
\Lam=\set{M_\e=M(u(\e))}{\abs\e<\e_0}
\end{equation}
will form a foliation, if we can show that
\begin{equation}
\dot u\ne0.
\end{equation}
Differentiating the equations in \re{5.70c} with respect to $\e$ and evaluating the result at $\e=0$ yields
\begin{equation}
A\dot u(0)=\dot\tau(0)\q\wed\q\int_M\h\dot u(0)=1
\end{equation}
and we deduce further
\begin{equation}
\dot\tau (0)\int_M\h=\spd{A\dot u(0)}\h=\spd{\dot u (0)}{A\h}=0
\end{equation}
and thus
\begin{equation}
\dot\tau(0)=0\q\wed\q\dot u(0)=\h>0,
\end{equation}
if $\h$ is normalized such that $\spd\h\h=1$, i.e., we have $\dot u(\e)>0$, if $\e_0$ is chosen small enough.
\ep

\br\lar{5.7.1}
Let $M$ be a stable solution of
\begin{equation}
\fv GM=0
\end{equation}
as in the preceding theorem, but not strictly stable and let $M_\e$ be a foliation of a neighbourhood of $M$ such that
\begin{equation}
\fv G{M_\e}=\tau(\e)\qq\A\,\abs\e<\e_0.
\end{equation}
If $M$ is the limit hypersurface of a curvature flow as in \rt{5.3c}, then
\begin{equation}
\tau(\e)>0\qq\A\, 0<\e<\e_0,
\end{equation}
if the flow hypersurfaces $M(t)$ converge to $M$ from above, which is tantamount to
\begin{equation}\lae{5.79c}
\F(F)-\tilde f\ge0,
\end{equation}
or we have
\begin{equation}
\tau(\e)<0\qq\A\,-\e_0<\e<0,
\end{equation}
if
\begin{equation}
\F(F)-\tilde f\le 0,
\end{equation}
in which case the flow hypersurfaces converge to $M$ from below.

The direction \cq{above} is defined by the region the normal $\s\nu$ of $M$ points to.
\er

\bp
Let us assume that the flow hypersurfaces satisfy \re{5.79c} and fix $0<\e<\e_0$. We may also suppose that the initial hypersurface $M(0)$ doesn't intersect the tubular neighbourhood of $M$ which is being foliated by $M_\e$.
Now, fix $0<\e<\e_0$, then there must be a first $t>0$ such that $M(t)$ touches $M_\e$ from above which yields, in view of the maximum principle,
\begin{equation}
\fv G{M_\e}=\tau(\e)>0,
\end{equation}
since $\tau(\e)\le 0$ would imply $\tau(\e)=0$ and $M_\e=M(t)$, \cf \cite[Theorem 2.7.9]{cg:cp}, i.e., $M(t)$ would be a stationary solution, which is impossible as we have proved in \rl{5.4c}.
\ep

Finally, let us show that the symmetric polynomials $H_k$, $1\le k\le n$, are of class $(D)$, if the ambient space has constant curvature.

\bl\lal{5.8.2}
Let $N$ be a semi-Riemannian space of constant curvature, then the symmetric polynomials $F=H_k$, $1\le k\le n$, are of class $(D)$. In case $k=2$ it suffices to assume $N$ Einstein.
\el

\bp
We shall prove the result by induction on $k$. First we note that the cones of definition $\C_k\su\R[n]$ of the $H_k$ form an ordered chain
\begin{equation}
\C_k\su\C_{k-1} \qq\A\, 1<k\le n,
\end{equation}
\cf \cite{garding}, so that a hypersurface admissible for $H_k$ is also admissible for $H_{k-1}$.

For $k=1$ we have
\begin{equation}
F^{ij}=g^{ij}
\end{equation}
and the result is obviously valid for arbitrary $N$.

Thus let us assume that the result is already proved for $1\le k<n$. Set $F=H_{k+1}$, $\hat F=H_k$ and let $M$ be an admissible hypersurface for $F$ with principal curvatures $\ka_i$.

From the definition of the $H_k$'s we immediately deduce 
\begin{equation}
\hat F=\pde F{\ka_i}+\ka_i\pde{\hat F}{\ka_i}
\end{equation}
for fixed $i$, no summation over $i$, or equivalently,
\begin{equation}
\hat Fg^{ij}=F^{ij}+\hat F^{jm}h_m^i,
\end{equation}
notice that the last term is a symmetric tensor, since for any symmetric curvature function $F$ $F^{ij}$ and $h_{ij}$ commute, \cf \cite[Lemma 2.1.9]{cg:cp}. Thus there holds
\begin{equation}
F^{ij}=\hat Fg^{ij}-\hat F^{jm}h_m^i
\end{equation}
and we deduce, using the induction hypothesis,
\begin{equation}
\begin{aligned}
F^{ij}_{;j}&=\hat F^i-\hat F^{jm}h^i_{m;j}=\hat F^i-\hat F^{jm}h_{mj;}^{\hp{mj;}i}\\
&=\hat F^i-\hat F^i=0,
\end{aligned}
\end{equation}
where we applied the Codazzi equations at one point.

If $F=H_2$, then
\begin{equation}
F^{ij}=Hg^{ij}-h^{ij}
\end{equation}
and the assumption $N$ Einstein suffices to conclude that $F^{ij}$ is divergence free.
\ep

\section{Existence results}

From now on we shall assume that ambient manifold $N$ is  Lorentzian, or more precisely, that it is smooth, globally hyperbolic with a  compact, connected Cauchy hypersurface. Then there exists a smooth future oriented time function $x^0$ such that the metric in $N$ can be expressed in Gaussian coordinates $(x^\al)$ as 
\begin{equation}\lae{6.1.1}
d\bar s^2=e^{2\psi}\{-(dx^0)^2+\s_{ij}dx^idx^j\},
\end{equation}
where $x^0$ is the time function and the $(x^i)$ are local coordinates for 
\begin{equation}
\so=\{x^0=0\}.
\end{equation}
$\so$ is then also a compact, connected Cauchy hypersurface. For a proof of the splitting result see \cite[Theorem 1.1]{sanchez:splitting}, and for the fact that all Cauchy hypersurfaces are diffeomorphic and hence $\so$ is also compact and connected, see \cite[Lemma 2.2]{Bernal:2003jb}.

One advantage of working in globally hyperbolic spacetimes with a compact Cauchy hypersurface is that all compact, connected  spacelike $C^m$-hypersurfaces $M$ can be written as graphs over $\so$.

\bl\lal{6.1.1}
Let $N$ be as above and $M\su N$ a connected, spacelike hypersurface of class $C^m$, $1\le m$, then $M$ can be written as a graph over $\so$
\begin{equation}
M=\graph \fv u{\so}
\end{equation}
with $u\in C^m(\so)$.
\el

We proved this lemma under the additional hypothesis that $M$ is achronal, \cite[Proposition 2.5]{cg:indiana},  however, this assumption is unnecessary as has been shown in \cite[Theorem 1.1]{heiko:diplom}.

We are looking at the curvature flow \fre{1.3.5} and want to prove that it converges to a stationary solution hypersurface, if certain assumptions are satisfied.

The existence proof consists of four steps:
\bir
\item
Existence on a maximal time interval $[0,T^*)$.

\item
Proof that the flow stays in a compact subset.

\item
Uniform a priori estimates in an appropriate function space, e.g., $C^{4,\al}(\so)$ or $C^{\un}(\so)$, which, together with (ii), would imply $T^*=\un$.

\item
Conclusion that the flow---or at least a subsequence of the flow hypersur\-faces---converges if $t$ tends to infinity.
\ei

The existence on a maximal time interval is always guaranteed, if the data are sufficiently regular, since the problem is parabolic. If the flow hypersurfaces can be written as graphs in a Gaussian coordinate system, as will always be the case in a globally hyperbolic spacetime with a compact Cauchy hypersurface in view of \rl{6.1.1}, the conditions are better than in the general case:
\bt\lat{6.2.2}
Let $4\le m\in\N$ and $0<\al<1$, and assume the semi-Riemannian space $N$ to be of class $C^{m+2,\al}$. Let the strictly monotone curvature function $F$, the functions $f$ and $\F$ be of class $C^{m,\al}$ and let $M_0\in C^{m+2,\al}$ be an admissible compact, spacelike, connected, orientable\footnote{Recall that oriented simply means there exists a continuous normal, which will always be the case in a globally hyperbolic spacetime.} hypersurface. Then the curvature flow \fre{1.3.5} with initial hypersurface $M_0$ exists in  a maximal time interval $[0,T^*)$, $0<T^*\le\un$, where in case that the flow hypersurfaces cannot be expressed as graphs they are supposed to be smooth, i.e, the conditions should be valid for arbitrary $4\le m\in\N$ in this case.
\et

A proof can be found in \cite[Theorem 2.5.19, Lemma 2.6.1]{cg:cp}.

\cvm
The second step, that the flow stays in a compact set, can only be achieved by barrier assumptions, \cf \rd{1.6.7}. Thus, let $\Om\su N$ be open and precompact such that $\pa\Om$ has exactly two components
\begin{equation}\lae{6.4.1}
\pa\Om=M_1\uud M_2
\end{equation}
where $M_1$ is a lower barrier for the pair $(F,f)$ and $M_2$ an upper barrier. Moreover, $M_1$ has to lie in the past of $M_2$
\begin{equation}\lae{6.5.1}
M_1\su I^-(M_2),
\end{equation}
\cf \cite[Remark 2.7.8]{cg:cp}.

Then the flow hypersurfaces will always stay inside $\bar\Om$, if the initial hypersurface $M_0$ satisfies $M_0\su \Om$, \cite[Theorem 2.7.9]{cg:cp}. This result is also valid if $M_0$ coincides with one the barriers, since then the velocity $(\F-\tilde f)$ has a weak sign and the flow moves into $\Om$ for small $t$, if it moves at all, and the arguments of the proof are applicable.

In Lorentzian manifolds the existence of barriers is associated with the presence of past and future singularities. In globally hyperbolic spacetimes, when $N$ is topologically a product
\begin{equation}
N=I\times\so,
\end{equation}
where $I=(a,b)$, singularities can only occur, when the endpoints of the interval are approached. A singularity, if one exists, is called a \tit{crushing singularity}, if the sectional curvatures become unbounded, i.e.,
\begin{equation}
\riema\al\bet\ga\de \bar R^{\al\bet\ga\de}\ra\un
\end{equation}
and such a singularity should provide a future \resp past barrier for the mean curvature function $H$.

\bd\lad{3.6.1.2}
Let $N$ be a globally hyperbolic spacetime with compact Cauchy hypersurface $\so$
so that $N$ can be written as a topological product $N=I\times \so$ and its
metric expressed as
\begin{equation}\lae{3.6.0.4}
d\bar s^2=e^{2\psi}(-(dx^0)^2+\s_{ij}(x^0,x)dx^idx^j).
\end{equation}
Here, $x^0$ is a globally defined future directed time function and $(x^i)$ are local
coordinates for $\so$.
$N$ is said to have a \tit{future} \resp \tit{past mean
curvature barrier}, if there are sequences $M_k^+$ \resp $M_k^-$ of closed,
spacelike, admissible  hypersurfaces such that
\begin{equation}\lae{6.9.1}
\lim_{k\ra\un} \fv H{M_k^+}=\un \q\tup{\resp}\q \lim_{k\ra\un} \fv H{M_k^-}=-\un
\end{equation}
and
\begin{equation}\lae{6.10.1}
\limsup \inf_{M_k^+}x^0> x^0(p)\qq\A\,p\in N
\end{equation}
\resp
\begin{equation}\lae{6.11.1}
\liminf \sup_{M_k^-}x^0< x^0(p)\qq\A\,p\in N,
\end{equation}
\ed

If one stipulates that the principal curvatures of the $M_k^+$ \resp $M_k^-$ tend to plus \resp minus infinity, then these hypersurfaces could also serve as barriers for other curvature functions. The past barriers would most certainly be non-admissible for any curvature function except $H$. 

\br\lar{6.4.2}
Notice that the assumptions \re{6.9.1} alone already implies \re{6.10.1} \resp \re{6.11.1}, if either
\begin{equation}\lae{6.12.1}
\limsup \inf_{M_k^+}x^0>a
\end{equation}
\resp
\begin{equation}
\liminf \sup_{M_k^-}x^0 <b
\end{equation}
where $(a,b)=x^0(N)$, or, if 
\begin{equation}\lae{6.14.1}
\bar R_{\al\bet}\nu^\al\nu^\bet\ge -\Lam\qq\A\, \spd\nu\nu=-1.
\end{equation}
where $\Lam\ge 0$. 
\er

\bp
It suffices to prove that the relation \re{6.10.1} is automatically satisfied under the assumptions \re{6.12.1} or \re{6.14.1}  by switching the light cone and replacing $x^0$ by $-x^0$ in case of the past barrier. 

Fix $k$, and let 
\begin{equation}
\tau_k=\inf_{M_k}x^0,
\end{equation}
 then the coordinate slice 
 \begin{equation}
M_{\tau_k}=\{x^0=\tau_k\}
\end{equation}
touches $M_k$ from below in a point $p_k\in M_k$ where $\tau_k=x^0(p_k)$ and the maximum principle yields that in that point
\begin{equation}
\fv H{M_{\tau_k}}\ge \fv H{M_k},
\end{equation}
hence, if $k$ tends to infinity the points $(p_k)$ cannot stay in a compact subset, i.e., 
\begin{equation}\lae{6.18.2}
\limsup x^0(p_k)\ra b
\end{equation}
or
\begin{equation}\lae{6.19.1}
\limsup x^0(p_k)\ra a.
\end{equation}

We shall show that only \re{6.18.2} can be valid. The relation \re{6.19.1} evidently contradicts \re{6.12.1}.

In case the assumption \re{6.14.1}  is valid, we consider a fixed coordinate slice $M_0=\{x^0=\const\}$, then all hypersurfaces $M_k$ satisfying
\begin{equation}
\fv H{M_0}<\inf_{M_k}H\q\wed\q \sqrt{n\Lam}<\inf_{M_k}H
\end{equation}
have to lie in the future of $M_0$, \cf \cite[Lemma 4.7.1]{cg:cp}, hence the result.
\ep

A future mean curvature barrier certainly represents a singularity, at least if $N$ 
satisfies the condition
\begin{equation}
\bar R_{\al\bet}\nu^\al\nu^\bet\ge -\Lam \qq\A\,\spd\nu\nu=-1
\end{equation}
where $\Lam\ge 0$, because of the future timelike incompleteness, which is proved in \cite{galloway:cft}, and is a generalization of Hawking's earlier result for $\Lam=0$, \cite{he:book}.  But these
singularities need not be crushing, \cf \cite[Section 2]{cg:imcf} for a counterexample.

\cvm
The uniform a priori estimates for the flow hypersurfaces are the hardest part in any existence proof. When the flow hypersurfaces can be written as graphs it suffices to prove $C^1$ and $C^2$ estimates, namely, the induced metric
\begin{equation}
g_{ij}(t,\xi)=\spd{x_i}{x_j}
\end{equation}
where $x=x(t,\xi)$ is a local embedding of the flow, should stay uniformly positive definite, i.e., there should exist positive constants $c_i$, $1\le i\le 2$, such that
\begin{equation}
c_1 g_{ij}(0,\xi)\le g_{ij}(t,\xi)\le c_2 g_{ij}(0,\xi),
\end{equation}
or equivalently, that the quantity
\begin{equation}
\tilde v=\spd{\h}\nu,
\end{equation}
where $\nu$ is the past directed normal of $M(t)$ and $\h$ the vector field
\begin{equation}
\h=(\h_\al)=e^\psi(-1,0,\ldots,0),
\end{equation}
is uniformly bounded, which is achieved with the help of the parabolic equation \fre{1.4.4.3}, if it is possible at all.

However, in some special situations $C^1$-estimates are automatically satisfied, \cf \rt{1.6.11} at the end of this section.

For the $C^2$-estimates  the principle curvatures $\ka_i$ of the flow hypersurfaces have to stay in a compact set in the cone of definition $\C$ of $F$, e.g., if $F$ is the Gaussian curvature, then $\C=\C_+$ and one has to prove that there are positive constants $k_i$, $i=1,2$ such that
\begin{equation}
k_1\le \ka_i\le k_2\qq\A\, 1\le i\le n
\end{equation}
uniformly in the cylinder $[0,T^*)\times M_0$, where $M_0$ is  any manifold that can serve as a base manifold for the embedding $x=x(t,\xi)$.

The parabolic equations that are used for these curvature estimates are \fre{1.4.3.13}, usually for an upper estimate, and \fre{1.4.3.10} for the lower estimate. Indeed, suppose that the flow starts at the upper barrier, then
\begin{equation}\lae{6.18.1}
F\ge f
\end{equation}
at $t=0$ and this estimate remains valid throughout the evolution because of the parabolic maximum principle, use \re{1.4.3.13}. Then, if upper estimates for the $\ka_i$ have been derived and if $f>0$ uniformly, then we conclude from \re{6.18.1} that the $\ka_i$ stay in a compact set inside the open cone $\C$, since 
\begin{equation}
\fv F{\pa\C}=0.
\end{equation}

To obtain higher order estimates we are going to exploit the fact that the flow hypersurfaces are graphs over $\so$ in an essential way, namely, we look at the associated scalar flow equation \fre{1.4.21} satisfied by $u$. This equation is a nonlinear uniformly parabolic equation, where the operator $\F(F)$ is also concave in $h_{ij}$, or equivalently, convex in $u_{ij}$, i.e., the $C^{2,\al}$-estimates of Krylov and Safonov, \cite[Chapter 5.5]{nk} or see \cite[Chapter 10.6]{oliver:pde} for a very clear and readable presentation, are applicable, yielding uniform estimates for the standard parabolic H\"older semi-norm
\begin{equation}\lae{6.29.2}
[D^2 u]_{\bet,\bar Q_T}
\end{equation}
for some $0<\bet\le\al$ in the cylinder
\begin{equation}
Q=[0,T)\times \so,
\end{equation}
independent of $0<T<T^*$, which in turn will lead to $H^{m+2+\al,\frac{m+2+\al}2}(\bar Q_T)$ estimates, \cf \cite[Theorem 2.5.9, Remark 2.6.2]{cg:cp}.  

$H^{m+2+\al,\frac{m+2+\al}2}(\bar Q_T)$ is a parabolic H\"older space, \cf \cite[p. 7]{lus} for the original definition and \cite[Note 2.5.4]{cg:cp} in the present context.

The estimate \re{6.29.2} combined with the uniform $C^2$-norm leads to uniform
 $C^{2,\bet}(\so)$-estimates independent of $T$.

These estimates imply that $T^*=\un$.

\cvm
Thus, it remains to prove that $u(t,\cdot)$ converges in $C^{m+2}(\so)$ to a stationary solution $\tilde u$, which is then also of class $C^{m+2,\al}(\so)$ in view of the Schauder theory.

Because of the preceding a priori estimates $u(t,\cdot)$ is precompact in $C^2(\so)$. Moreover, we deduce from the scalar flow equation \fre{1.4.21} that $\dot u$ has a sign, i.e., the $u(t,\cdot)$ converge monotonely in $C^0(\so)$ to $\tilde u$ and therefore also in $C^2(\so)$.

To prove that $\graph \tilde u$ is a solution, we again look at \re{1.4.21} and integrate it with respect to $t$ to obtain for fixed $x\in\so$
\begin{equation}
\abs{\tilde u(x)-u(t,x)}=\int_t^\un e^{-\psi}v\abs{\F-\tilde f},
\end{equation}
where we used that $(\F-\tilde f)$ has a sign, hence $(\F-\tilde f)(t,x)$ has to vanish when $t$ tends to infinity, at least for a subsequence, but this suffices to conclude that $\graph \tilde u$ is a stationary solution and
\begin{equation}
\lim_{t\ra\un}(\F-\tilde f)=0.
\end{equation}

Using the convergence of $u$ to $\tilde u$ in $C^2$, we can then prove:
\bt\lat{6.5.5}
The functions $u(t,\cdot)$ converge in $C^{m+2}(\so)$ to $\tilde u$, if the data satisfy the assumptions in \rt{6.2.2}, since we have
\begin{equation}\lae{6.33.2}
u\in H^{m+2+\bet,\frac{m+2+\bet}2}(\bar Q),
\end{equation}
where $Q=Q_\un$. 
\et

\bp
Out of convention let us write $\al$ instead of $\bet$ knowing that $\al$ is the H\"older exponent in \re{6.29.2}.

We shall reduce the Schauder estimates to the standard Schauder estimates in $\R[n]$ for the heat equation with a right-hand side by using the already established results \re{6.29.2} and 
\begin{equation}
u(t,\cdot)\underset{C^2(\so)}\ra \tilde u\in C^{m+2,\al}(\so).
\end{equation}

Let $(U_k)$ be a finite open covering of $\so$ such that each $U_k$ is contained in a coordinate chart and
\begin{equation}\lae{6.35.2}
\diam U_k<\rho,
\end{equation}
$\rho$ small, $\rho$ will be specified in the proof, and let $(\h_k)$ be a subordinate finite partition of unity of class $C^{m+2,\al}$.

Since
\begin{equation}
u\in H^{m+2+\al,\frac{m+2+\al}2}(\bar Q_T)
\end{equation}
for any finite $T$, \cf \cite[Lemma 2.6.1]{cg:cp}, and hence 
\begin{equation}
u(t,\cdot)\in C^{m+2,\al}(\so)\qq\A\,0\le t<\un
\end{equation}
we shall choose $u_0=u(t_0,\cdot)$ as initial value for some large $t_0$ such that
\begin{equation}\lae{6.38.2}
\abs{a^{ij}(t,\cdot)-\tilde a^{ij}}_{0,\so}<\e/2\qq\A\,t\ge t_0,
\end{equation}
where
\begin{equation}
a^{ij}=v^2\dot\F F^{ij}
\end{equation}
and  $\tilde a^{ij}$ is defined correspondingly for $\tilde M=\graph \tilde u$.

However, making a variable transformation we shall always assume that $t_0=0$ and $u_0=u(0,\cdot)$.

We shall prove \re{6.33.2} successively.

\cvm
(i) Let us first show that
\begin{equation}
D_xu\in H^{2+\al,\frac{2+\al}2}(\bar Q).
\end{equation}
This will be achieved, if we show that for an arbitrary $\xi\in C^{m+1,\al}(T^{1,0}(\so))$
\begin{equation}
\f=D_\xi u\in H^{2+\al,\frac{2+\al}2}(\bar Q),
\end{equation}
\cf \cite[Remark 2.5.11]{cg:cp}.

Differentiating the scalar flow equation \fre{1.4.21} with respect to $\xi$ we obtain
\begin{equation}\lae{6.42.2}
\dot \f-a^{ij}\f_{ij}+b^i\f_i+c\f=f,
\end{equation}
where of course the symbol $f$ has a different meaning then in \re{1.4.21}.

Later we want to apply the Schauder estimates for solutions of the heat flow equation with right-hand side. In order to use elementary potential estimates we have to cut off $\f$ near the origin $t=0$ by considering
\begin{equation}
\tilde\f=\f \theta,
\end{equation}
where $\theta=\theta(t)$ is smooth satisfying
\begin{equation}
\theta(t)=
\begin{cases}
1,&t>1,\\
0,&t\le\frac12.
\end{cases}
\end{equation}
This modification doesn't cause any problems, since we already have a priori estimates for finite $t$, and we are only concerned about the range $1\le t<\un$. $\tilde\f$ satisfies the same equation as $\f$ only the right-hand side has the additional summand $w\dot\theta$.

Let $\h=\h_k$ be one of the members of the partition of unity and set
\begin{equation}
w=\tilde\f\h,
\end{equation}
then $w$ satisfies a similar equation with slightly different right-hand side
\begin{equation}\lae{6.44.2}
\dot w-a^{ij}w_{ij}+b^iw_i+cw=\tilde f
\end{equation}
but we shall have this in mind when applying the estimates.

The $w(t,\cdot)$ have compact support in one of the $U_k$'s, hence we can replace the covariant derivatives of $w$ by ordinary partial derivatives without changing the structure of the equation and the properties of the right-hand side, which still only depends linearly on $\f$ and $D\f$.

We want to apply the well-known estimates for the ordinary heat flow equation
\begin{equation}
\begin{aligned}
\dot w-\D w&=\hat f
\end{aligned}
\end{equation}
where $w$ is defined  in $\R[]\times\R[n]$.

To reduce the problem to this special form, we pick an arbitrary $x_0\in U_k$, set $z_0=(0,x_0)$, $z=(t,x)$ and consider instead of \re{6.44.2}
\begin{equation}
\begin{aligned}
\dot w-a^{ij}(z_0)w_{ij}&=\hat f\\
&=[a^{ij}(z)-a^{ij}(z_0)]w_{ij}-b^iw_i-cw+\tilde f,
\end{aligned}
\end{equation}
where we emphasize that the difference
\begin{equation}\lae{6.47.2}
\abs{a^{ij}(z)-a^{ij}(z_0)}
\end{equation}
can be made smaller than any given $\e>0$ by choosing $\rho=\rho(\e)$ in \re{6.35.2} and $t_0=t_0(\e)$ in \re{6.38.2} accordingly. Notice also that this equation can be extended into $\R[]\times \R[n]$, since all functions have support in $\{t\ge \frac12\}$.

Let $0<T<\un$ be arbitrary, then all terms belong to the required function spaces in $\bar Q_T$ and there holds
\begin{equation}
[w]_{2+\al,Q_T}\le c[\hat f]_{\al,Q_T},
\end{equation}
where $c=c(n,\al)$. The brackets indicate the standard unweighted parabolic semi-norms, \cf \cite[Definition 2.5.2]{cg:cp}, which are identical to those defined in \cite[p. 7]{lus}, but there the brackets are replaced by kets. 

Thus, we conclude
\begin{equation}
\begin{aligned}
\hp{a}&[w]_{2+\al,Q_T}\le c \sup_{U_k\times (0,T)}\abs{a^{ij}(z)-a^{ij}(z_0)}[D^2w]_{\al,Q_T}+c[f]_{\al,Q_T} \\
&+c_1\{[D^2u]_{\al,Q_+}+[Du]_{\al,Q_T}+[u]_{\al,Q_T}+\abs{w}_{0,Q_T}+\abs{D^2w}_{0,Q_T}\},
\end{aligned}
\end{equation}
where $c_1$ is independent of $T$, but dependent on $\h_k$. Here we also used the fact that the lower order coefficients and $\f, D\f$ are uniformly bounded.

Choosing now $\e>0$ so small that
\begin{equation}
c\msp \e<\tfrac12
\end{equation}
and $\rho$, $t_0$ accordingly such the difference in \re{6.47.2} is smaller than $\e$, we deduce
\begin{equation}
\begin{aligned}
\hp{a}&[w]_{2+\al,Q_T}\le 2c[f]_{\al,Q_T} \\
& +2c_1\{[D^2u]_{\al,Q_+}+[Du]_{\al,Q_T}+[u]_{\al,Q_T}+\abs{w}_{0,Q_T}+\abs{D^2w}_{0,Q_T}\}.
\end{aligned}
\end{equation}

Summing over the partition of unity and noting that $\xi$ is arbitrary we see that in the preceding inequality we can replace $w$ by $Du$ everywhere  resulting in the estimate
\begin{equation}
\begin{aligned}
&\hp{a}[Du]_{2+\al,Q_T}\le 
c_1[f]_{\al,Q_T}\\
&\qq +c_1\{[D^2u]_{\al,Q_T}+[Du]_{\al,Q_T}+[u]_{\al,Q_T}+\abs{Du}_{0,Q_T}+\abs{D^3u}_{0,Q_T}\},
\end{aligned}
\end{equation}where $c_1$ is a new constant still independent of $T$.

Now the only critical terms on the right-hand side are $\abs{D^3u}_{0,Q_T}$, which can be estimated by \re{6.55.2}, and the H\"older semi-norms with respect to $t$
\begin{equation}
[Du]_{\frac\al 2,t,Q_T}+[u]_{\frac\al 2,t,Q_T}.
\end{equation}
The second  one is taken care of by the boundedness of $\dot u$, see \fre{1.4.21}, while the first one is estimated with the help of equation \re{6.42.2} revealing
\begin{equation}
\abs{D\dot u}\le c\{\sup_{[0,T]}\abs{u}_{3,\so}+\abs{f}_{0,Q_T}\},
\end{equation}
since for fixed but arbitrary $t$ we have
\begin{equation}\lae{6.55.2}
\abs{u}_{3,\so}\le \e[D^3u]_{\al,\so}+c_\e\abs{u}_{0,\so},
\end{equation}
where $c_\e$ is independent of $t$.

Hence we conclude
\begin{equation}
\abs{Du}_{2+\al,Q_T}\le\const
\end{equation}
uniformly in $T$.

\cvm
(ii) Repeating these estimates successively for $2\le l\le m$ we obtain uniform estimates for
\begin{equation}
\sum_{l=2}^m[D^l_xu]_{2+\al,Q_T},
\end{equation}
which, when combined with the uniform $C^2$-estimates, yields
\begin{equation}\lae{6.58.2}
\abs{u(t,\cdot)}_{m+2,\al,\so}\le \const
\end{equation}
uniformly in $0\le t<\un$.

Looking at the equation \re{1.4.21} we then deduce
\begin{equation}\lae{6.59.2}
\abs{\dot u(t,\cdot)}_{m,\al,\so}\le \const
\end{equation}
uniformly in $t$.

\cvm
(iii) To obtain the estimates for $D^r_tu$ up to the order
\begin{equation}
[\tfrac{m+2+\al}2]
\end{equation}
we differentiate the scalar curvature  equation with respect to $t$ as often as necessary and also with respect to the mixed derivatives $D^r_tD^s_x$ to estimate
\begin{equation}
\sum_{1\le 2r+s<m+2+\al}D^r_tD^s_x u
\end{equation}
using \re{6.58.2}, \re{6.59.2} and the results from the prior differentiations. 

Combined with the estimates for the heat equation in $\R[]\times\R[n]$ these estimates will also yield the necessary a priori estimates for the H\"older semi-norms in $\bar Q$, where again the smallness of \re{6.47.2} has to be used repeatedly.
\ep

\br
The preceding regularity result is also valid in Riemannian manifolds, if the flow hypersurfaces can be written as graphs in a Gaussian coordinate system. In fact the proof is unaware of the nature of the ambient space.
\er

With the method described above the following existence results have been proved in globally hyperbolic spacetimes with a compact Cauchy hypersurface. $\Om\su N$ is always a precompact domain the boundary of which is decomposed as in \re{6.4.1} and \re{6.5.1} into an upper and lower barrier for the pair $(F,f)$. We also apply the stability results from \rs{5c} and the just proved regularity of the convergence and formulate the theorems accordingly. 

By convergence of the flow in $C^{m+2}$ we mean convergence of the leaves $M(t)=\graph u(t,\cdot)$ in this norm.

\bt
Let $M_1$, $M_2$ be lower \resp upper barriers for the pair $(H,f)$, where $f\in C^{m,\al}(\bar\Om)$ and the $M_i$ are of class $C^{m+2,\al}$, $4\le m$, $0<\al<1$, then the curvature flow
\begin{equation}
\begin{aligned}
\dot x&=(H-f)\nu\\
x(0)&=x_0,
\end{aligned}
\end{equation}
where $x_0$ is an embedding of the initial hypersurface $M_0=M_2$ exists for all time and converges in $C^{m+2}$ to a stable solution $M$ of class $C^{m+2,\al}$ of the equation
\begin{equation}
\fv HM=f,
\end{equation}
provided the initial hypersurface is not already a solution.
\et

The existence result was proved in \cite[Theorem 2.2]{cg:mz}, see also \cite[Theorem 4.2.1]{cg:cp} and the remarks following the theorem. Notice that $f$ isn't supposed to satisfy any sign condition.

For spacetimes that satisfy the timelike convergence condition and for functions $f$ with special structural conditions existence results via a mean curvature flow were first proved in \cite{eh1}.

The Gaussian curvature or the curvature functions $F$ belonging to the larger class $(K^*)$, see \cite{cg:indiana} for a definition, require that the admissible hypersurfaces are strictly convex. 

Moreover, proving a priori estimates for the second fundamental form of a hypersurface $M$ in general semi-Riemannian manifolds, when the curvature function is not the mean curvature, or does not behave similar to it,  requires that a strictly convex function $\chi$ is defined in a neighbourhood of the hypersurface, see \frl{2.8.1} where sufficient assumptions are stated which imply the existence of strictly convex  functions.

Furthermore, when we consider curvature functions of class $(K^*)$, notice that the Gaussian curvature belongs to that class, then the right-hand side $f$ can be defined in $T(\bar\Om)$ instead of $\bar\Om$, i.e., in a local trivialization of the tangent bundle $f$ can be expressed as
\begin{equation}
f=f(x,\nu)\q\wed\q\nu\in T_x(N).
\end{equation}
We shall formulate the existence results with this more general assumption, though of course any stability claim only makes sense for $f=f(x)$.

\bt
Let $F\in C^{m,\al}(\C_+)$, $4\le m$, $0<\al<1$, be a curvature function of class $(K^*)$, let $0<f\in C^{m,\al}(T(\bar\Om))$, and let $M_1$, $M_2$ be lower \resp upper barriers for $(F,f)$ of class $C^{m+2,\al}$. Then the curvature flow
\begin{equation}
\begin{aligned}
\dot x&=(\F-\tilde f)\nu\\
x(0)&=x_0
\end{aligned}
\end{equation}
where $\F(r)=\log r$ and $x_0$ is an embedding of $M_0=M_2$, exists for all time and converges in $C^{m+2}$ to a stationary solution $M\in C^{m+2,\al}$ of the equation
\begin{equation}
\fv FM=f
\end{equation}
provided the initial hypersurface $M_2$ is not already a stationary solution and there exists a strictly convex function $\chi\in C^2(\bar\Om)$.

When $f=f(x)$ and $F$ is of class $(D)$, then $M$ is stable.
\et

The theorem was proved in \cite{cg:indiana} when $f$ is only defined in $\bar\Om$ and in the general case in \cite[Theorem 4.1.1]{cg:cp}.

\cvm
When $F=H_2$ is the scalar curvature operator, then the requirement that $f$ is defined in the tangent bundle and not merely in $N$ is a necessity, if the scalar curvature is to be prescribed. To prove existence results in this case, $f$ has to satisfy some natural structural conditions, namely,
\begin{align}
0<c_1&\le f(x,\n)\qq\tup{if}\q\spd\n\n=-1,\lae{0.3}\\[\cma]
\nnorm{f_\bet(x,\n)}&\le c_2 (1+\nnorm\n^2),\lae{0.4}\\
\intertext{and}
\nnorm{f_{\n^\bet}(x,\n)}&\le c_3 (1+\nnorm\n),\lae{0.5}
\end{align}
\nd for all $x\in\bar\Om$ and all past directed timelike vectors $\n\in T_x(\Om)$,
where $\nnorm{\cdot}$ is a Riemannian reference metric.

Applying a curvature flow to obtain stationary solutions requires to approximate $F$ and $f$ by functions $F_\e$ and $f_k$ and to use these functions for the flow. The $F_\e$ are the $\e$-regularizations of $F$, which we already discussed before, \cf \fre{5.31.1}. Let us also write $\tilde F$ instead of $F_\e$ as before. 

The functions $f_k$ have the property that $\nnorm{f_{k\bet}}$ only grows linearly in $\nnorm\nu$ and $\nnorm{f_{k\n^\bet}(x,\n)}$ is bounded. To simplify the presentation we shall therefore assume that $f$ satisfies
\begin{equation}\lae{6.33.1}
\nnorm{f_\bet(x,\n)}\le c_2 (1+\nnorm\n),
\end{equation}
\begin{equation}\lae{6.34.1}
\nnorm{f_{\n^\bet}(x,\n)}\le c_3,
\end{equation}
and also 
\begin{equation}\lae{6.35.1}
0<c_1\le f(x,\n)\qq\A\,\nu\in T_x(N),\;\spd\nu\nu<0,
\end{equation}
although the last assumption is only a minor point that can easily be dealt with, see \cite[Remark 2.6]{cg:scalar}, and \cite[Section 7 and 8]{cg:scalar} for the other approximations of $f$.

The barriers $M_i$, $i=1,2$, for $(F,f)$ satisfy the barrier condition of course only weakly, i.e., no strict inequalities; however, because of the $\e$-regularization we need strict inequalities, so that the $M_i$'s are also barriers for $(\tilde F,f)$, if $\e$ is small. In \cite[Remark 2.4 and Lemma 2.5]{cg:scalar} it is shown that strict inequalities for the barriers may be assumed without loss of generality.

Now, we can formulate the existence result for the scalar curvature operator $F=H_2$ under these provisions.
\bt\lat{6.6.1}
Let $f\in C^{m,\al}(T(\bar\Om))$, $4\le m$, $0<\al<1$, satisfy the conditions \re{6.33.1}, \re{6.34.1} and \re{6.35.1}, and let $M_1$, $M_2$ be strict lower \resp upper barriers of class $C^{m+2,\al}$ for $(F,f)$. Let $\tilde F$ be the $\e$-regularization of the scalar curvature operator $F$, then the curvature flow for $\tilde F$
\begin{equation}
\begin{aligned}
\dot x&=(\F-\tilde f)\\
x(0)&=x_0
\end{aligned}
\end{equation}
where $\F(r)=r^\frac12$ and $x_0$ is an embedding of $M_0=M_2$, exists for all time and converges in $C^{m+2}$ to a stationary solution $M_\e\in C^{m+2,\al}$ of
\begin{equation}
\fv{\tilde F}{M_\e}=f
\end{equation}
provided there exists a strictly convex function $\chi\in C^2(\bar\Om)$ and $0<\e$ is small.

The $M_\e$ then converge in $C^{m+2}$ to a solution $M\in C^{m+2,\al}$ of
\begin{equation}
\fv FM=f.
\end{equation}
If $f=f(x)$ and $N$ Einstein, then $M$ is stable.  
\et

These statements, except for the stability and the convergence in $C^{m+2}$, are proved in \cite{cg:scalar}.

\br
Let us now discuss the pure mean curvature flow
\begin{equation}
\dot x=H\nu
\end{equation}
with initial spacelike hypersurface $M_0$ of class $C^{m+2,\al}$, $m\ge 4$ and $0<\al<1$. From the corresponding scalar curvature flow \fre{1.4.21} we immediately infer that the flow moves into the past of $M_0$, if
\begin{equation}\lae{6.40.1}
\fv H{M_0}\ge0
\end{equation}
and into its future, if
\begin{equation}
\fv H{M_0}\le0.
\end{equation}
Let us only consider the case \re{6.40.1} and also assume that $M_0$ is not maximal. From the a priori estimates in \cite[Section 3 and Section 4]{cg:mz} we then deduce that the flow remains smooth as long as it stays in a compact set of $N$, and if a compact, spacelike hypersurface $M_1$ of class $C^2$ satisfying
\begin{equation}\lae{6.42.1}
\fv H{M_1}\le0
\end{equation}
lies in the past of $M_0$, then the flow will exist for all time and converge in $C^{m+2}$ to a stable maximal hypersurface $M$, hence a neighbourhood of $M$ can be foliated by CMC hypersurfaces, where those in the future of $M$ have positive mean curvature, in view of \frr{5.7.1}.

Thus, the flow will converge if and only if such a hypersurface $M_1$ lies in the past of $M_0$.

Conversely, if there exists a compact, spacelike hypersurface $M_1$ in $N$ satisfying \re{6.42.1}, and there is no stable maximal hypersurface in its future, then this is a strong indication that $N$ has no future singularity, assuming that such a singularity would produce spacelike hypersurfaces with positive mean curvature.

An example of such a spacetime is the $(n+1)$-dimensional de Sitter space which is geodesically complete and has exactly one maximal hypersurface $M$ which is also totally geodesic but not stable, and the future \resp past of $M$ are foliated by coordinate slices with negative \resp positive mean curvature.
\er

To conclude this section let us show which spacelike hypersurfaces satisfy $C^1$-estimates automatically.

\bt\lat{1.6.11}
Let $M=\graph \fv u\so$ be a compact, spacelike hypersurface represented in a
Gaussian coordinate system with  unilateral bounded  principal curvatures,
e.g.,
\begin{equation}
\kappa_i\ge \kappa_0\q\A\,i.
\end{equation}
Then, the quantity $\tilde v=\frac{1}{\sql}$ can be estimated by
\begin{equation}
\tilde v\le c(\abs u,\so,\sigma_{ij},\psi,\kappa_0),
\end{equation}
where we assumed that in the Gaussian coordinate system the
ambient metric has the form as in \re{6.1.1}.
\et
\bp
We suppose as usual that the Gaussian coordinate system is future oriented, and that the
second fundamental form is evaluated with respect to the past directed normal.
We observe that
\begin{equation}\lae{1.6.87}
\norm{Du}^2=g^{ij}u_iu_j=e^{-2\psi}\frac{\abs{Du}^2}{v^2}\raise 2pt \hbox{,}
\end{equation}
hence, it is equivalent to find an a priori estimate for $\norm{Du}$.

Let $\lambda$ be a real parameter to be specified later, and set
\begin{equation}
w=\tfrac{1}{2}\log\norm{Du}^2+\lambda u.
\end{equation}
We may regard $w$ as being defined on $\so$; thus, there is $x_0\in\so$ such that
\begin{equation}
w(x_0)=\sup_\so w,
\end{equation}
and we conclude
\begin{equation}
0=w_i=\frac{1}{\norm{Du}^2}\,u_{ij}u^j+\lambda u_i
\end{equation}
in $x_0$, where the covariant derivatives are taken with respect to the induced
metric
$g_{ij}$, and the indices are also raised with respect to that metric.

Expressing the second fundamental form of a graph with the help of the Hessian of the function
\begin{equation}\lae{2.16}
e^{-\psi}v^{-1}h_{ij}=-u_{ij}-\cha 000\mspace{1mu}u_iu_j-\cha 0i0
\mspace{1mu}u_j-\cha 0j0\mspace{1mu}u_i-\cha ij0.
\end{equation}
 we deduce further
\begin{equation}\lae{4.16}
\begin{aligned}
\lambda\norm{Du}^4&=-u_{ij}u^iu^j\\
&= e^{-\psi}\tilde vh_{ij}u^iu^j+\cha 000\msp \norm{Du}^4\\
&\hp{=}\msp[2]+2\cha 0j0\msp u^j\norm{Du}^2+\cha ij0\msp u^iu^j.
\end{aligned}
\end{equation}
Now, there holds
\begin{equation}
u^i=g^{ij}u_j=e^{-2\psi}\sigma^{ij}u_jv^{-2},
\end{equation}
and by assumption,
\begin{equation}
h_{ij}u^iu^j\ge \kappa_0\msp\norm{Du}^2,
\end{equation}
i.e., the critical terms on the right-hand side of \re{4.16} are of fourth order in
$\norm{Du}$ with bounded coefficients, and we conclude that $\norm{Du}$ can't be
too large in $x_0$ if we choose $\lambda$ such that
\begin{equation}
\lambda\le -c\msp\nnorm{\cha \alpha\beta 0}-1
\end{equation}
with a suitable constant $c$; $w$, or equivalently, $\norm{Du}$ is therefore
uniformly bounded from above.
\ep

Especially for convex graphs over $\so$ the term $\tilde v$ is uniformly bounded as long as
they stay in a compact set.  

\section{The inverse mean curvature flow}
Let us now consider the inverse mean curvature flow (IMCF)
\begin{equation}\lae{7.1}
\dot x=-H^{-1}\nu
\end{equation}
with initial hypersurface $M_0$  in a globally hyperbolic spacetime $N$ with compact Cauchy hypersurface $\so$.

$N$ is supposed to satisfy the timelike convergence condition
\begin{equation}
\bar R_{\al\bet}\nu^\al\nu^\bet\ge0\qq\A\,\spd\nu\nu=-1.
\end{equation}
Spacetimes with compact Cauchy hypersurface that satisfy the timelike convergence condition are also called \tit{cosmological spacetimes}, a terminology due to Bartnik.

In such spacetimes the inverse mean curvature flow will be smooth as long as it stays in a compact set, and, if $\fv H{M_0}>0$ and if the flow exists for all time, it will necessarily run into the future singularity, since the mean curvature of the flow hypersurfaces will become unbounded and the flow will run into the future of $M_0$. Hence the claim follows from \frr{6.4.2}.

However, it might be that the flow will run into the singularity in finite time. To exclude this behaviour we introduced in \cite{cg:imcf} the so-called \tit{strong volume
decay condition}, \cf \rd{5.1.2}. A strong volume decay condition is both
necessary and sufficient in order that the IMCF exists for all time.

\bt\lat{7.1}
Let $N$ be a cosmological spacetime with compact Cauchy hypersurface $\so$ and
with a  future mean curvature barrier. Let $M_0$ be a closed, connected, spacelike hypersurface
with positive mean curvature and assume furthermore that $N$ satisfies a future
volume decay condition. Then the IMCF \re{7.1} with initial hypersurface $M_0$ exists
for all time and provides a foliation of the future \inds{$D^+(M_0)$} of $M_0$.

The evolution parameter $t$ can be chosen as a new time function. The flow
hypersurfaces $M(t)$ are the slices $\{t=\const\}$ and their volume satisfies
\begin{equation}
\abs{M(t)}=\abs{M_0} e^{-t}.
\end{equation}

Defining a new time function $\tau$ by choosing
\begin{equation}
\tau=1-e^{-\frac1n t}
\end{equation}
we obtain  $0\le \tau <1$,
\begin{equation}
\abs{M(\tau)}=\abs{M_0} (1-\tau)^n,
\end{equation}
and  the future singularity corresponds to $\tau=1$.

Moreover, the length $L(\ga)$ of any future directed curve $\ga$ starting from
$M(\tau)$ is bounded from above by
\begin{equation}
L(\ga)\le c (1-\tau),
\end{equation}
where $c=c(n, M_0)$. Thus, the expression $1-\tau$ can be looked at as the radius
of the slices $\{\tau=\const\}$ as well as a measure of the remaining life span of the
spacetime.
\et

Next we shall define the \ind{strong volume decay condition}.

\bd\lad{5.1.2}
Suppose there exists a time function $x^0$ such that the future end of $N$ is
determined by $\{\tau_0\le x^0<b\}$ and the coordinate slices
$M_\tau=\{x^0=\tau\}$ have positive mean curvature with respect to the past
directed normal for $\tau_0\le\tau<b$. In addition the volume $\abs{M_\tau}$
should satisfy
\begin{equation}\lae{5.1.19}
\lim_{\tau\ra b}\abs{M_\tau}=0.
\end{equation}

A decay like that is normally associated with a future singularity and we simply call it
\tit{\ind{volume decay}}. If $(g_{ij})$ is the induced metric of $M_\tau$ and
$g=\det(g_{ij})$, then we have
\begin{equation}\lae{5.1.20}
\log g(\tau_0,x)-\log g(\tau,x)=\int_{\tau_0}^\tau 2 e^\psi \bar H(s,x)\q\A\,x\in \so,
\end{equation}
where $\bar H(\tau,x)$ is the mean curvature of $M_\tau$ in $(\tau,x)$. This relation can be easily derived from the relation \fre{1.3.6} and \frr{1.3.5}.  A detailed proof
is given in \cite{cg:volume}.

In view of \re{5.1.19} the left-hand side of this equation tends to infinity if $\tau$
approaches $b$ for \aev $x\in \so$, i.e.,
\begin{equation}
\lim_{\tau\ra b}\int_{\tau_0}^\tau e^\psi \bar H(s,x)=\un\q \tup{for
\aev}\; x\in\so.
\end{equation}

Assume now, there exists a continuous, positive function $\f=\f(\tau)$ such that
\begin{equation}\lae{5.1.22}
e^\psi \bar H(\tau,x)\ge \f(\tau)\qq\A\, (\tau, x)\in (\tau_0,b)\times \so,
\end{equation}
where
\begin{equation}\lae{5.1.23}
\int_{\tau_0}^b \f(\tau)=\un,
\end{equation}
then we say that the future of $N$ satisfies a \tit{strong volume decay condition}.
\ed

\br
(i) By approximation we may assume that the function $\f$ above is
smooth.

(ii) A similar definition holds for the past of $N$ by simply reversing the time
direction. Notice that in this case the mean curvature of the coordinate slices has to
be negative.
\er

\bl\lal{5.1.4}
Suppose that the future of $N$ satisfies a strong volume decay condition, then there
exist a time function $\tilde x^0=\tilde x^0(x^0)$, where $x^0$ is the time function
in the strong volume decay condition, such that the mean curvature $\bar H$ of the
slices $\tilde x^0=\const$ satisfies the estimate
\begin{equation}\lae{5.1.24}
e^{\tilde\psi}\bar H\ge 1.
\end{equation}
The factor $e^{\tilde\psi}$ is now the conformal factor in the representation
\begin{equation}\lae{5.1.25}
d\bar s^2=e^{2\tilde\psi}(-(d\tilde x^0)^2+\s_{ij}dx^idx^j).
\end{equation}

The range of $\tilde x^0$ is equal to the interval $[0,\un)$, i.e., the singularity
corresponds to $\tilde x^0=\un$.
\el

A proof is given in \cite[Lemma 1.4]{cg:imcf}.

\br
\rt{7.1} can be generalized to spacetimes satisfying
\begin{equation}
\bar R_{\al\bet}\nu^\al\nu^\bet\ge-\Lam\qq\A\,\spd\nu\nu=-1
\end{equation}
with a constant $\Lam\ge 0$, if the mean curvature of the initial hypersurface $M_0$ is sufficiently large
\begin{equation}
\fv H{M_0}>\sqrt{n\Lam},
\end{equation}
\cf \cite{heiko:diplom}. In that thesis it is also shown that the future mean curvature barrier assumption can be dropped, i.e., the strong volume decay condition is sufficient to prove that the IMCF exists for all time and provides a foliation of the future of $M_0$. Hence, the strong volume decay condition already implies the existence of a future mean curvature barrier, since the leaves of the IMCF define  such a barrier.
\er

\section{The IMCF in ARW spaces}
In the present section we consider spacetimes $N$ satisfying some structural
conditions, which are still fairly general, and prove convergence results for the leaves
of the IMCF.

Moreover, we define a new spacetime $\hat N$ by switching the light
cone and using reflection to define a new time function, such that the two
spacetimes $N$ and $\hat N$ can be pasted together to yield a smooth manifold
having a metric singularity, which, when viewed from the region $N$ is a big crunch,
and when viewed from $\hat N$ is a big bang.

The inverse mean curvature flows in $N$ \resp $\hat N$ correspond to each other
via reflection. Furthermore, the properly rescaled flow in $N$ has a natural smooth
extension of class $C^3$ across the singularity into $\hat N$. With respect to this
natural diffeomorphism we speak of a transition from big crunch to big bang. 

\bd\lad{6.0.1}
A globally hyperbolic spacetime $N$, $\dim N=n+1$, is said to be \tit{\ind{asymptotically
Robertson-Walker}} (ARW) with respect to the future, if a future end of $N$, $N_+$,
can be written as a product $N_+=[a,b)\times \so$, where $\so$ is a
Riemannian space, and there exists a future directed time function $\tau=x^0$ such
that the metric in $N_+$ can be written as
\begin{equation}\lae{6.1.19}
d\breve s^2=e^{2\tilde\psi}\{-{(dx^0})^2+\s_{ij}(x^0,x)dx^idx^j\},
\end{equation}
where  $\so$ corresponds to $x^0=a$, $\tilde\psi$ is of the form
\begin{equation}
\tilde\psi(x^0,x)=f(x^0)+\psi(x^0,x),
\end{equation}
and we assume that there exists a positive constant $c_0$ and a smooth
Riemannian metric $\bar\s_{ij}$ on $\so$ such that
\begin{equation}
\lim_{\tau\ra b}e^\psi=c_0\q\wed\q \lim_{\tau\ra b}\s_{ij}(\tau,x)=\bar\s_{ij}(x),
\end{equation}
and
\begin{equation}
\lim_{\tau\ra b}f(\tau)=-\un.
\end{equation}

Without loss of generality we shall assume $c_0=1$. Then $N$ is ARW with
respect to the future, if the metric is close to the Robertson-Walker metric
\begin{equation}\lae{6.0.5}
d\bar s^2=e^{2f}\{-{dx^0}^2+\bar\s_{ij}(x)dx^idx^j\}
\end{equation}
near the singularity $\tau =b$. By \tit{close} we mean that the derivatives of arbitrary order with respect to space and time of the
conformal metric $e^{-2f}\breve g_{\al\bet}$ in \re{6.1.19} should converge  to the
corresponding derivatives of the conformal limit metric in \re{6.0.5} when $x^0$ tends
to $b$. We emphasize that in our terminology \ind{Robertson-Walker metric} does not
imply that
$(\bar\s_{ij})$ is a metric of constant curvature, it is only the spatial metric of a
warped product.

We assume, furthermore, that $f$ satisfies the following five conditions
\begin{equation}
-f'>0,
\end{equation}
there exists $\om\in\R[]$ such that
\begin{equation}\lae{6.0.7}
n+\om-2>0\q\wed\q \lim_{\tau\ra b}\abs{f'}^2e^{(n+\om-2)f}=m>0.
\end{equation}
Set $\tilde\ga =\frac12(n+\om-2)$, then there exists the limit
\begin{equation}\lae{6.0.8}
\lim_{\tau\ra b}(f''+\tilde\ga \abs{f'}^2)
\end{equation}
and
\begin{equation}\lae{6.0.9}
\abs{D^m_\tau(f''+\tilde\ga \abs{f'}^2)}\le c_m \abs{f'}^m\qq
\A\, m\ge 1,
\end{equation}
as well as
\begin{equation}\lae{6.0.10}
\abs{D_\tau^mf}\le c_m \abs{f'}^m\qq\A\, m\ge 1.
\end{equation}

If $\so$ is compact, then we call $N$ a \tit{normalized} ARW spacetime,\index{normalized ARW spacetime} if
\begin{equation}
\int_{\so}\sqrt{\det{\bar\s_{ij}}}=\abs{S^n}. 
\end{equation}
\ed

\br
(i) If these assumptions are satisfied, then  the range of $\tau$ is
finite, hence, we may---and shall---assume \wlogc that $b=0$, i.e.,
\begin{equation}
a<\tau<0.
\end{equation}

(ii) Any ARW spacetime with compact $\so$ can be normalized as one easily checks. For normalized ARW
spaces the constant
$m$ in \re{6.0.7} is defined uniquely and can be identified
with the mass of $N$, \cf \cite{cg:mass}.

(iii) In view of the assumptions on $f$ the mean curvature of the coordinate slices
$M_\tau=\{x^0=\tau\}$ tends to $\un$, if $\tau$ goes to zero.

(iv) ARW spaces with compact $\so$ satisfy a strong volume decay condition, \cf
\frd{5.1.2}.

(v) Similarly one can define $N$ to be ARW with respect to the past. In this case the
singularity would lie in the past, correspond to $\tau=0$, and the mean curvature
of the coordinate slices would tend to $-\un$. 
\er

We assume that $N$ satisfies the timelike convergence condition and that $\so$ is compact. Consider the future
end $N_+$ of $N$ and let $M_0\su N_+$ be a spacelike hypersurface with positive
mean curvature $\fv {\breve H}{M_0}>0$ with respect to the past directed normal
vector $\breve\nu$---it will become apparent in  a moment why we use the symbols $\breve H$
and $\breve\nu$ and not the usual ones $H$ and $\nu$. Then, as we have proved in
the preceding section, the inverse mean curvature flow
\begin{equation}
\dot x=-\breve H^{-1}\breve\nu
\end{equation}
with initial hypersurface $M_0$ exists for all time, is smooth, and runs straight
into the future singularity.

If we express the flow hypersurfaces $M(t)$  as graphs over $\so$
\begin{equation}
M(t)=\graph u(t,\cdot),
\end{equation}
then we have proved in \cite{cg:arw}

\bt\lat{6.0.3}
\tup{(i)} Let $N$ satisfy the above assumptions, then the range of the time function
$x^0$ is finite, i.e., we may assume that $b=0$. Set
\begin{equation}\lae{8.15}
\tilde u=ue^{\ga t},
\end{equation}
where $\ga=\tfrac1n\tilde\ga$, then there are positive constants $c_1, c_2$ such
that
\begin{equation}
-c_2\le \tilde u\le -c_1<0,
\end{equation}
and $\tilde u$ converges in $C^\un(\so)$ to a smooth function, if $t$ goes to
infinity. We shall also denote the limit function  by $\tilde u$.

\tup{(ii)} Let $\breve g_{ij}$ be the induced metric of the leaves $M(t)$, then the
rescaled metric
\begin{equation}
e^{\frac2n t} \breve g_{ij}
\end{equation}
converges in $C^\un(\so)$ to
\begin{equation}
(\tilde\ga m)^\frac1{\tilde\ga}(-\tilde u)^\frac2{\tilde\ga}\bar\s_{ij}.
\end{equation}

\tup{(iii)} The leaves $M(t)$ get more umbilical, if $t$ tends to infinity, namely, there
holds
\begin{equation}
\breve H^{-1}\abs{\breve h^j_i-\tfrac1n \breve H\de^j_i}\le c\msp e^{-2\ga t}.
\end{equation}
In case $n+\om-4>0$, we even get a better estimate
\begin{equation}
\abs{\breve h^j_i-\tfrac1n \breve H\de^j_i}\le c\msp e^{-\frac1{2n}(n+\om-4) t}.
\end{equation}
\et

To prove the convergence results for the inverse mean curvature flow, we
consider the flow hypersurfaces to be embedded in $N$ equipped with the conformal
metric
\begin{equation}\lae{6.2.1}
d\bar s^2=-(dx^0)^2+ \s_{ij}(x^0,x)dx^idx^j.
\end{equation}

Though, formally, we have a different ambient space we still denote it by the same
symbol $N$ and distinguish only the metrics $\breve g_{\al\bet}$ and $\bar
g_{\al\bet}$
\begin{equation}
\breve g_{\al\bet}=e^{2\tilde\psi}\bar g_{\al\bet}
\end{equation}
and the corresponding geometric quantities of the hypersurfaces $\breve h_{ij},
\breve g_{ij}, \breve\nu$ \resp $h_{ij}, g_{ij}, \nu$, etc., i.e., the standard  notations  now apply  to the case when $N$ is equipped with the metric in \re{6.2.1}.

The second fundamental forms $\breve h^j_i$ and $h^j_i$ are related by
\begin{equation}
e^{\tilde\psi}\breve h^j_i=h^j_i+\tilde\psi_\al\nu^\al\de^j_i
\end{equation}
and, if we define $F$ by
\begin{equation}
F=e^{\tilde\psi} \breve H,
\end{equation}
then
\begin{equation}
F=H-n\tilde vf'+n\psi_\al\nu^\al,
\end{equation}
where
\begin{equation}
\tilde v=v^{-1},
\end{equation}
and the evolution equation can be written as
\begin{equation}\lae{6.2.7}
\dot x=-F^{-1}\nu,
\end{equation}
since
\begin{equation}
\breve\nu=e^{-\tilde\psi}\nu.
\end{equation}

The flow exists for all time and is smooth, due to the results in the preceding section.

Next, we want to show how the metric, the second fundamental form, and the
normal vector of the hypersurfaces $M(t)$ evolve by adapting the general evolution equations in \frs{1.3} to the present situation. 
\bl\lal{6.2.1}
The metric, the normal vector, and the second fundamental form of $M(t)$
satisfy the evolution equations
\begin{equation}\lae{6.2.9}
\dot g_{ij}=-2 F^{-1}h_{ij},
\end{equation}
\begin{equation}\lae{6.2.10}
\dot \n=\nabla_M(-F^{-1})=g^{ij}(-F^{-1})_i x_j,
\end{equation}
and
\begin{equation}\lae{6.2.11}
\dot h_i^j=(-F^{-1})_i^j+F^{-1} h_i^k h_k^j + F^{-1} \riema
\al\bet\ga\de\n^\al x_i^\bet \n^\ga x_k^\de g^{kj}
\end{equation}
\begin{equation}
\dot h_{ij}=(-F^{-1})_{ij}-F^{-1} h_i^k h_{kj}+ F^{-1} \riema
\al\bet\ga\de\n^\al x_i^\bet \n^\ga x_j^\de.
\end{equation}
\el

Since the initial hypersurface is a graph over $\so$, we can write
\begin{equation}
M(t)=\graph\fu{u(t)}S0\q \A\,t\in I,
\end{equation}
where $u$ is defined in the cylinder $\R[]_+\times \so$. We then deduce from
\re{6.2.7}, looking at the component $\al=0$, that $u$ satisfies a parabolic
equation of the form
\begin{equation}\lae{6.2.14}
\dot u=\frac{\tilde v}F,
\end{equation}
where  we emphasize that the time
derivative is a total derivative, i.e.
\begin{equation}
\dot u=\pde ut+u_i\dot x^i.
\end{equation}

Since the past directed normal can be expressed as
\begin{equation}
(\n^\al)=-e^{-\psi}v^{-1}(1,u^i),
\end{equation}
we conclude from \re{6.2.14}
\begin{equation}\lae{6.2.17}
\pde ut=\frac vF.
\end{equation}

For this new curvature flow the necessary decay estimates and convergence results can be proved, which in turn can be immediately translated to corresponding convergence results for the original IMCF.

\hinweis*{Transition from big crunch to big bang}
With the help of the convergence results  in \rt{6.0.3}, we can rescale the IMCF such that it can be extended past the singularity in a natural way.

We  define a new spacetime $\hat N$ by reflection and time reversal such that
the IMCF in the old spacetime transforms to an IMCF in the new one.

By switching the light cone we obtain a new spacetime $\hat N$. The flow equation
in $N$ is independent of the time orientation, and we can write it as
\begin{equation}
\dot x=-\breve H^{-1}\breve\nu=-(-\breve H)^{-1}(-\breve\nu)\equiv -\hat
H^{-1}\hat \nu,
\end{equation}
where the normal vector $\hat \nu=-\breve\nu$ is past directed in $\hat N$ and the
mean curvature $\hat H=-\breve H$ negative.

Introducing a new time function $\hat x^0=-x^0$ and formally new coordinates
$(\hat x^\al)$ by setting
\begin{equation}
\hat x^0=-x^0,\q\hat x^i=x^i,
\end{equation}
we define a spacetime $\hat N$ having the same metric as $N$---only expressed in
the new coordinate system---such that the flow equation has the form
\begin{equation}\lae{6.8.3}
\dot{\hat x}=-\hat H^{-1}\hat \nu,
\end{equation}
where $M(t)=\graph \hat u(t)$, $\hat u=-u$, and 
\begin{equation}
(\hat\nu^\al)=-\tilde ve^{-\tilde \psi}(1,\hat u^i)
\end{equation}
in the new coordinates, since
\begin{equation}
\hat\nu^0=-\breve \nu^0\pde{\hat x^0}{x^0}=\nu^0
\end{equation}
and
\begin{equation}
\hat\nu^i=-\breve\nu^i.
\end{equation}

The singularity in $\hat x^0=0$ is now a past singularity, and can be referred to as a
big bang singularity.

The union $N\uu\hat N$ is a smooth manifold, topologically a product
$(-a,a)\ti\so$---we are well aware that formally the singularity $\{0\}\ti\so$ is not
part of the union; equipped with the respective metrics and time orientation it is a
spacetime which has a (metric) singularity in
$x^0=0$. The time function
\begin{equation}\lae{6.8.7}
\hat x^0=
\begin{cases}
\hp{-}x^0, &\tup{in } N,\\
-x^0, &\tup{in } \hat N,
\end{cases}
\end{equation}
is smooth across the singularity and future directed.

$N\uu\hat N$ can be regarded as a \tit{\ind{cyclic universe}} with a contracting part
$N=\{\hat x^0<0\}$ and an expanding part $\hat N=\{\hat x^0>0\}$ which are
joined at the singularity $\{\hat x^0=0\}$.

It turns out that the inverse mean curvature flow, properly rescaled, defines a
natural $C^3$\nobreak- diffeomorphism across the singularity and with respect to
this diffeomorphism we speak of a transition from big crunch to big bang.

Using the time function in \re{6.8.7} the inverse mean curvature flows in $N$ and
$\hat N$ can be uniformly expressed in the form
\begin{equation}\lae{6.8.8}
\dot{\hat x}=-\hat H^{-1}\hat\nu,
\end{equation}
where \re{6.8.8} represents the original flow in $N$, if $\hat x^0<0$, and the flow in
\re{6.8.3}, if $\hat x^0>0$.

Let us now introduce a new flow parameter
\begin{equation}
s=
\begin{cases}
-\ga^{-1}e^{-\ga t},& \tup{for the flow in } N,\\
\hp{-}\ga ^{-1}e^{-\ga t},& \tup{for the flow in } \hat N,
\end{cases}
\end{equation}
and define the flow $y=y(s)$ by $y(s)=\hat x(t)$. $y=y(s,\x)$ is then defined in
$[-\ga^{-1},\ga^{-1}]\times \so$, smooth in $\{s\ne 0\}$, and satisfies the
evolution equation
\begin{equation}\lae{6.8.10}
y'\equiv \tfrac d{ds}y=
\begin{cases}
-\hat H^{-1}\hat\nu \msp e^{\ga t}, & s<0,\\
\hp{-}\hat H^{-1}\hat\nu \msp e^{\ga t},& s>0.
\end{cases}
\end{equation}

In \cite{cg:arw} we proved:
\bt\lat{6.8.1}
The flow $y=y(s,\x)$  is of class $C^3$ in $(-\ga^{-1},\ga^{-1})\times
\so$ and defines a natural diffeomorphism across the singularity. The flow parameter
$s$ can be used as a new time function.
\et

\br
The regularity result for the transition flow is optimal, i.e., given any $0<\al<1$, then there is an ARW space such that the transition flow is not of class $C^{3,\al}$, \cf \cite{cg:rw}.
\er

\br\lar{8.4}
Since ARW spaces have a future mean curvature barrier, a future end can be foliated by CMC hypersurfaces the mean curvature of which can be used as a new time function., see \cite{cg1} and \cite{cg:foliation2}. In \cite{cg:arwcmc} we study this foliation a bit more closely and prove that, when writing the CMC hypersurfaces as graphs $M_\tau=\graph \f(\tau,\cdot)$ in the special coordinate system of the ARW space, where $\tau$ is the mean curvature, of $M_\tau$ then
\begin{equation}
\tau (-\f)^{1+\tilde\ga^{-1}}\ra \const>0,
\end{equation}
notice that $\f<0$, and hence
\begin{equation}
\lim_{\tau\ra\un}\frac{\f(\tau,x)}{\f(\tau,y)}=1\qq\A\,x,y\in\so.
\end{equation}

Moreover, the new time function
\begin{equation}
s=-\tau^{-q},\q q=\frac{\tilde\ga}{1+\tilde\ga}
\end{equation}
can be extended to the mirror universe $\hat N$ by odd reflection as a function of class $C^3$ across the singularity with non-vanishing gradient.
\er


 
 \providecommand{\bysame}{\leavevmode\hbox to3em{\hrulefill}\thinspace}
\providecommand{\href}[2]{#2}



\end{document}